%% file: StabilityAIRGA.tex
\documentclass[final]{siamart1116}

\input{ex_shared}

\ifpdf
\hypersetup{
	pdftitle={\TheTitle},
	pdfauthor={\TheAuthors}
}
\fi

\begin{document}
	
	%
	%
	%
	%


	\maketitle
	
	\begin{abstract}
		{Simulation of large dynamical systems can be unmanageable due to high demands on computational resources. These large systems can be reduced into a smaller dimension by using Model Order Reduction (MOR) techniques. The reduced system has approximately the same characteristics as the original system but it requires significantly less computational effort in simulation. MOR can be done in many ways such as balanced truncation, Hankel approximations, and Krylov projection. Among these, the projection methods are quite popular, and hence, we focus on them.}			
		{{There exist many projection methods based MOR algorithms for reducing an extensive range of linear dynamical systems.} That is, non-parametric-parametric as well as first-order and second-order. 
			
		{Here, we focus on MOR of non-parametric second-order dynamical systems.} In these MOR algorithms, sequences of very large and sparse linear systems arise during the model reduction process. Solving such linear systems is the main computational bottleneck in efficient scaling of these MOR algorithms for reducing extremely large dynamical systems.}
		{Preconditioned iterative methods are often used for solving such linear systems.}
		 
		{These iterative methods introduce errors because they solve the linear systems up to a certain tolerance. {Hence, our focus} is to analyze the stability of 
		{the selected category of} MOR algorithms {(non-parametric second-order)} when using inexact linear solves. {Adaptive Iterative Rational Global Arnoldi (AIRGA) \cite{Bonin20161} is a popular MOR algorithm belonging to this category}.}
		{
			We prove that, under {four} mild  conditions, the AIRGA algorithm is backward stable with respect to the errors introduced by these inexact linear solves. {Our results easily extend to other MOR algorithms belonging to this category.} 
			{Our \textit{first condition} enforces the use of a Ritz-Galerkin based linear solver, where the
				residual of a linear system is made orthogonal to the corresponding Krylov subspace. 
				Our \textit{second condition} requires satisfying few extra orthogonalities. Since Conjugate Gradient (CG) is the most popular method based upon the Ritz-Galerkin theory, we use it. We show how to modify CG to achieve these extra orthogonalities.			
				Modifying CG with the suggested changes is non-trivial. Hence, we further demonstrate that using Recycling CG (RCG) helps us achieve these orthogonalities with no code changes. The extra cost of orthogonalizations is often offset by savings because of recycling.} 
				
				{Our \textit{third condition} involves existence and invertibility of a matrix mostly dependent upon the input dynamical system, with the norm of this matrix bounded by one.		
					 Our \textit{fourth and final condition} involves being able to compute a perturbation from the derived expression and bounding its norm by one as well. 
The last two conditions are easily satisfied by all our models.
					 }		
		}    
	\end{abstract}
	
	
	\section{Introduction}
	\label{intro}
	
{Dynamical systems arise in many areas of science and engineering. There are three factor that define a dynamical system; (i) linearity; (ii) parametrization; and (iii) order. Linear dynamical systems usually approximate the real-life phenomenas well, and hence, have been extensively studied. Thus, we focus on {\it linear} dynamical systems.}

{Whether a dynamical systems in parametrized or not; and the order of derivatives in the system are the other two characteristics defining a dynamical system.}
{In general}, a parameterized second-order dynamical system is usually of the form \cite{Feng2013}
\begin{align}\label{eq:Linear_param_system}
	\begin{split}
		M(p_1, p_2, \ldots, p_w)\ddot{x}(t) + D(p_1, p_2, \ldots, p_w)\dot{x}(t) + K(p_1, p_2, \ldots, p_w)x(t)  = Bu(t), \\ 
		y(t)  = C x(t),
	\end{split}
\end{align}
where $M(\cdot), D(\cdot), K(\cdot) \in \mathbb{R}^{n \times n}$, $B \in \mathbb{R}^{n \times m}$, $C \in \mathbb{R}^{q \times n}$ and $(p_1,  p_2,  \ldots,  p_w)$ are the parameters that are linearly embedded in the dynamical system matrices. Also, $x(t) \colon \mathbb{R} \rightarrow \mathbb{R}^n$ is the vector of all states,  $u(t) \colon \mathbb{R} \rightarrow \mathbb{R}^{m}$ and $y(t) \colon \mathbb{R} \rightarrow \mathbb{R}^{q}$ are the inputs and the outputs of the system, respectively. If $M(\cdot) = 0$, then above equation can be written in the form of the parametric first-order dynamical system as
\begin{align}\label{eq:first-Linear_param_system_1}
	\begin{split}
		D(p_1,  p_2,  \ldots,  p_w)\dot{x}(t) + K(p_1,  p_2,  \ldots,  p_w)x(t)  = Bu(t), \\ 
		y(t)  =  C x(t).
	\end{split}
\end{align}

If in \eqref{eq:Linear_param_system} and \eqref{eq:first-Linear_param_system_1} the system matrices are independent of the parameters, then they represent a non-parametric second-order and first-order dynamical system, respectively.

Simulation of large dynamical systems can be unmanageable due to high demands on computational resources. These large systems can be reduced into a smaller dimension by using Model Order Reduction (MOR) techniques \cite{Grimme1997,  Antoulas05,  Gugercin2008, Breiten2013}. The reduced system has approximately  the same characteristics as the original system but it requires significantly less computational effort in simulation. MOR can be done in many ways such as balanced truncation, Hankel approximations, and Krylov projection. Among these, the projection methods are quite popular, and hence, we focus on them.
Table \ref{tab:MOR-algo} summarizes most of the commonly used such algorithms.


\begin{table}[]
	\centering
	\footnotesize
	\caption{Linear MOR Algorithms based upon Projection.}
	\begin{tabular}{|c|c|c|c|}
		\hline
		S. No. & \textbf{Category} &  \textbf{Second-order} &  \textbf{First-order}\\ \hline
		\multirow{2}{*}{\textbf{1}} & \multirow{4}{*}{Parametric} & \cellcolor{gray!15}Cell 1\label{cell:par_1st} & \cellcolor{gray!15}Cell 2\label{cell:par_2nd} \\ \cline{3-4} 
		&  &  \begin{tabular}[c]{@{}c@{}}S-RPMOR \cite{Feng2013}, \\ IDPA \cite{SAADVANDI2014}, \\ S-PBTMR \cite{sonN2017} \end{tabular} &  \begin{tabular}[c]{@{}c@{}} IPMOR \cite{Baur2011}, \\ RPMOR \cite{Benner2014}, \\ PBTMR \cite{sonN2017} \end{tabular}\\ \hline
		\multirow{2}{*}{\textbf{2}} & \multirow{4}{*}{Non-Parametric} & \cellcolor{gray!15}Cell 3\label{cell:non_par_1st} & \cellcolor{gray!15}Cell 4\label{cell:non_par_2nd} \\ \cline{3-4} 
		&  & \begin{tabular}[c]{@{}c@{}}SOR-IRKA \cite{Wyatt2012}, \\  SO-IRKA \cite{Qiu2018},\\ SOSPDR \cite{BaiS05}, \\  AIRGA \cite{Bonin20161} \end{tabular} &  \begin{tabular}[c]{@{}c@{}} IRKA \cite{Gugercin2008}, \\ $(Sy)^{2}\text{IRKA}$ \cite{Breiten2013},\\  MIRIAm \cite{Gerstner2010} \end{tabular}\\ \hline
		
	\end{tabular}
	\label{tab:MOR-algo}
\end{table}

{In the mentioned MOR algorithms in Table \ref{tab:MOR-algo}},  sequences of very large and sparse linear systems arise during the model reduction process. Solving such linear systems is the main computational bottleneck in efficient scaling of these MOR algorithms for reducing extremely large dynamical systems, {which we discuss next.}
%

\subsection{Iterative Methods and Preconditioners}
\label{sec:prec_iter}

Direct methods, which are based upon different matrix factorizations, are commonly used for solving linear systems of equations \cite{Benzi2002}. Standard direct methods scale badly in-terms of the number of operations and the memory used (with respect to the increase in the size of the linear systems; as here). They typically perform dense linear algebra operations, and hence, are not an efficient choice when the linear system matrices are sparse (as here as well). 

An alternative to this is to use sparse direct methods. These methods solve this scaling problem to a great extent such that linear systems of fairly large size could be efficiently solved\footnote{{Often, they work well for linear systems arising from certain problem classes  \cite{Benzi2002,davis2006,Alon2013}, for example, discretization of PDEs in two dimensions.}}. However, sparse direct methods also become prohibitively expensive for extremely large sizes (hundreds of millions of equations to billions of equations). 

In such cases, using iterative methods are usually the only viable option, which scale well both in time and memory. Although iterative methods are not as robust or reliable as direct methods, they are still preferred {when} scaling is a bigger issue. {This is the case with MOR algorithms, and hence, we use them here.}

{\color{black}Krylov subspace based methods are very popular class of iterative methods \cite{Saad2003}, {which we focus on}. 
	{If} $
	Ax = b
	$ 
	{is the} linear system {to be solved}, with $A \in \mathbb{R}^{n \times n}, \ b \in \mathbb{R}^n$, $x_0$ the initial solution and $r_0$ (where $r_0 = b-Ax_0$) the initial residual, {then these methods}
	find the solution in $\mathbb{K}_{\mathcal{k}}(A, \ r_0) = span\{r_0, \ Ar_0, \ A^2r_0, \ \ldots, \ A^{{\mathcal{k}}-1}r_0\}$, where $\mathbb{K}_{\mathcal{k}}(\cdot, \ \cdot)$ represents the Krylov subspace.
	
	Often iterative methods are slow or fail to converge, and hence, preconditioning is used to accelerate them. We expect that the preconditioned iterative solves would find a solution in less amount of time as compared to the unpreconditioned ones. For most of the input dynamical systems, the Krylov subspace methods fail to converge. Hence, we use a preconditioner. The goal is to find a preconditioner that is cheap to compute as well as apply.

	If $P$ is a non-singular matrix that approximates the inverse of $A$, then the preconditioned system becomes $AP\tilde{x} = b$ with $x = P\tilde{x}$. 
	{\color{black}This is termed as right preconditioning. Similarly, left preconditioning can also be performed, where the preconditioner is present on 
		the left side of the matrix \cite{Benzi2002}\footnote{\color{black}If the preconditioner is present on both the sides of the coefficient matrix, then it is called split/  center preconditioning.}.
		If the linear system coefficient matrices are SPD, then both the types of preconditioning give the same results \cite{Benzi2002}. 
		
		For our MOR
		algorithms under-consideration, the linear system coefficient matrices do not have any special structure. Hence, both these
		types of preconditioning work differently.		
		In our experiments, we {mostly} use right preconditioning because it is
		fairly common \cite{Chow1998,Alexander2007}. However, to demonstrate that our techniques are independent of the type of preconditioning, for {some models}, we experiment with left preconditioning in the side as well.	
	}
	
	
	Preconditioned iterative methods are not exact because they solve linear systems upto a certain tolerance. {This raises} the question that if preconditioned iterative methods are used inside the MOR algorithms, then are these algorithms stable with respect to the error introduced by these methods. 
	Hence, {our focus is to} investigate the stability of MOR algorithms (with respect to use of iterative methods). This is briefly elaborated upon in Section \ref{sec:stab_intro}.
	
	
	\subsection{Stability Analysis of MOR Algorithms}\label{sec:stab_intro}
	As {mentioned} earlier, we investigate the stability of MOR algorithms with respect to use of preconditioned iterative methods. 
	This kind of 
	analysis was first proposed in \cite{Beattie20122916}, where a popular MOR algorithm for {\it linear} non-parametric first-order dynamical systems was analyzed (corresponding to Cell 4 of Table \ref{tab:MOR-algo}). 

	{In this paper, we focus on stability analysis of MOR algorithm belonging to Cell 3 of Table 1, that is,  for the non-parametric second-order. 
		Specifically, we focus on the Adaptive Iterative Rational Global Arnoldi (AIRGA) \cite{Bonin20161} algorithm with our results easily carried over to the other algorithms of this category. The stability of MOR algorithms belonging to the Cells 2 and 1 of Table \ref{tab:MOR-algo}, i.e. parametric first-order and second-order, respectively has been recently {dealth} in the dissertation of the second author \cite{navneet_thesis}.  
		}
	
	{An extended stability analyses for commonly used MOR algorithms for {\it bilinear} dynamical systems {(different parametrizations and orders)} have been done in \cite{CHOUDHARY201856} and \cite{Choudhary2019}. As earlier, our focus is on {\it linear} dynamical systems and not {\it bilinear}.}
	In the current context, it is important to highlight the difference between our {track of} stability analyses and the one done in \cite{LU2016} as well. The authors in \cite{LU2016} first showed that the SOAR algorithm {(for MOR of non-parametric second-order dynamical systems)} is unstable with respect to the machine precision errors (and not inexact solves of iterative methods, which is our focus). Then, they proposed a Two-level orthogonal Arnoldi (TOAR) algorithm that cures this instability of SOAR (we propose recycling variants of the underlying iterative methods for achieving stability).
	
	%

	The main contributions of this paper are as follows: We discuss stability of the AIRGA algorithm with respect to these inexact linear solves in Section \ref{sec:stab_anal_nonparametric}. 
	{
		{
			This {paper} has a very unique contribution that has not yet been looked at by any other past work. In \cite{Beattie20122916}, the authors mention that their stability analysis for non-parametric first-order dynamical systems can be easily carried to non-parametric second-order systems. Besides the fact that the authors do not perform this analysis in-details, they also do not focus on how to satisfy the arising stability conditions, which we do.}	
		{In all our analyses, we show that satisfying the stability conditions requires {\it changing the underlying linear solvers}, and that too in an efficient way so as to not incur any extra cost. {These aspects are discussed in Section \ref{sec:Sat_Back_stab_cond}}.}
	In Section \ref{sec:accuracy}, we derive the expression for accuracy of the reduced system, in-terms of the conditioning of the linear system as well as the residuals of the linear solves.}
	Numerical experiments, which support our preconditioned iterative solver
	theory are given in Section \ref{sec:Numerical_results}. Finally, we give conclusions and future directions in Section \ref{sec:Conc_future}.
	
	
	{The following sets of notations are used in this {paper}:} 
	\begin{itemize}
		\item $\mathbb{R}$ denotes the set of real numbers.
		\item $\mathbb{C}$ denotes the set of complex numbers.
		\item $\mathbb{N}$ denotes the set of natural numbers.
		\item $\|\cdot\|_{f}$ denotes the Frobenius norm.
		\item $\|\cdot\|$  denotes the Euclidean norm for vectors and the induced spectral norm for matrices.
		\item $I$ is the identity matrix.
	\end{itemize}

\section{Stability Analysis of AIRGA}\label{sec:stab_anal_nonparametric}
{The AIRGA algorithm is used for reducing non-parametric second-order} dynamical systems with {\it proportional damping}. {These systems have the form}

\begin{align}\label{eq:linear-non-para}
	\begin{split}
		M \ddot{x}(t) + D \dot{x}(t) + K x(t) & = B u(t), \\
		y(t) & =  Cx(t) =  \mathscr{C}^T x(t),
	\end{split}  
\end{align}        	   	 
where $M,\ D,\ K \in  \mathbb{R}^{n \times n}, \ B \in  \mathbb{R}^{n \times m}, \ \mathscr{C} \in  \mathbb{R}^{n \times q},$ and $D = \alpha M + \beta K$. Here, $\alpha, \ \beta$ are some scalar values. 
Let $V \in \mathbb{R}^{n \times r}$ and its columns span a $r$-dimension subspace ($r \ll n$).	
In principle, the Ritz-Galerkin projection method, {as used by AIRGA,} involves the steps below.

\begin{itemize}
	\item Approximating the reduced state vector $\hat{x}(t)$ using ${V}$ as $ x(t) \approx V \hat{x}(t)$ leads to
	\begin{align*}
		\begin{split}
			M V\ddot{\hat{x}}(t) + D V \dot{\hat{x}}(t)+KV \hat{x}(t) - B u(t) & = r(t),\\
			\hat{y}(t) & = \ \mathscr{C}^T V \hat{x}(t),
		\end{split}
	\end{align*}
	where $r(t)$ is the residual after projection. 
	\item Enforcing the residual $r(t)$ to be orthogonal to $V$ or $V^T r(t) = 0$ 
	leads to the reduced system  given as follows:
	{\begin{align}\label{TD_r1}
			\begin{split}
				& V^{T}\left(M V\ddot{\hat{x}}(t) + D V \dot{\hat{x}}(t)+KV \hat{x}(t) - B 
				u(t)\right)= 0,\\
				& \hat{y}(t)= \ \mathscr{C}^TV \hat{x}(t).
			\end{split}
		\end{align}}
		or 
		\begin{align}\label{eq:red_sys}
			\begin{split}
				\hat{M} \ddot{\hat{x}}(t) + \hat{D} \dot{\hat{x}}(t)+ \hat{K} \hat{x}(t) - \hat{B}	u(t) & = 0,\\
				\hat{y}(t) & =  \hat{\mathscr{C}}^T  \hat{x}(t),
			\end{split}
		\end{align}	
	\end{itemize}
	where 
	\begin{align}\label{eq:non_para_red}
		\hat{M} = V^{T}MV,\ \hat{D} = V^{T}DV,\ \hat{K}=V^{T}KV,\ \hat{B} = V^{T}B, \ \text{and} \ \hat{\mathscr{C}}^T = \mathscr{C}^TV. 
	\end{align}
	
	To compute this projection matrix $V$, AIRGA matches the moments of the original system transfer function and the reduced system transfer function.
	
	The transfer function of (\ref{eq:linear-non-para}) is given by 
	\begin{align*}
		H(s)=\mathscr{C}^T\left(s^{2}M+sD+K\right)^{-1}B= \mathscr{C}^TX(s),
	\end{align*}
	{where $X(s) = \left(s^{2}M+sD+K\right)^{-1}B.$}
	The power series expansion of $X\left(s\right)$ around an expansion point $s_{0} \in \mathbb{R}$ is given by (see, e.g.,~\cite{Wang2002})
	\begin{align}\label{eq:LT_eq1}
		X(s)=\sum\limits_{j=0}^{\infty} X^{(j)}(s_{0})\left(s-s_{0}\right)^{j},
	\end{align}
	where,
	\begin{align}\label{eq:Moment_orig}
		\begin{split}
			X^{(0)}\left(s_{0}\right)=& \ \left(s_{0}^{2}M+s_{0}D+K\right)^{-1}B,  \\ 
			X^{(1)}\left(s_{0}\right)=& \ \left(s_{0}^{2}M+s_{0}D+K\right)^{-1}\left(-\left(2s_{0}M+D\right)\right)X^{(0)}\left(s_{0}\right), \qquad \mathrm{and} \\ 
			X^{(j)}\left(s_{0}\right)= & \ \left(s_{0}^{2}M+s_{0}D+K\right)^{-1}\left[-\left(2s_{0}M+D\right)X^{(j-1)}\left(s_{0}\right)  -MX^{(j-2)}(s_{0})\right], 
		\end{split}
	\end{align}
	for $j = 2, 3, \ldots$. Here, $X^{(j)}\left(s_{0}\right)$ is called the $j^{th}$-order system moment at $s_{0}$.
	\par Similarly, the transfer function of the reduced system (\ref{eq:red_sys}) is given by 
	\begin{align*}
		\hat{H}(s)= \hat{\mathscr{C}}^T\hat{X}(s),
	\end{align*}
	where $\hat{X}(s)= \left(s^{2}\hat{M}+s\hat{D}+\hat{K} \right)^{-1}\hat{B}.$
	The power series expansion of $\hat{X}\left(s\right)$ around an expansion point $s_{0} \in \mathbb{R}$ is given by
	\begin{align}\label{eq:LT_r}
		\hat{X}(s)=\sum\limits_{j=0}^{\infty} \hat{X}^{(j)}(s_{0})\left(s-s_{0}\right)^{j}.
	\end{align}
	The $j^{th}$-order system moment $\hat{X}^{(j)}(s_{0})$ is defined analogously to $X^{(j)}(s_{0})$ in (\ref{eq:Moment_orig}). 
	
	\par The goal of moment-matching approach is to find a reduced system such that the first few moments of (\ref{eq:LT_eq1}) and (\ref{eq:LT_r}) are matched, that is, $X^{(j)}(s_0) = \hat X^{(j)}(s_0)$ for $j = 0, 1, 2, \ldots, t $ for some $t.$ This can be achieved by the observation below. With 
	\begin{align*}
		{P}_{1} = & - \left(s_{0}^{2}M+s_{0}D+K\right)^{-1}\left(2s_{0}M+ D\right), \\
		{P}_{2} = &- \left(s_{0}^{2}M+s_{0}D+K\right)^{-1}M, \\
		{Q} = & \ \left(s_{0}^{2}M+s_{0}D+K\right)^{-1}B,
	\end{align*}
	we have  from (\ref{eq:Moment_orig}) 
	\begin{align*}
		X^{(0)}\left(s_{0}\right) = & \  {Q},  \\
		X^{(1)}\left(s_{0}\right) = & \ {P}_{1}X^{(0)}\left(s_{0}\right), \qquad \mathrm{and}\\
		X^{(j)}\left(s_{0}\right)= & \ {P}_{1}  X^{(j-1)}\left(s_{0}\right) +  {P}_{2}X^{(j-2)}\left(s_{0}\right)
	\end{align*}
	for $j \ge 2.$ As already observed in \cite{BaiS05}, these moments are just the blocks of the second-order Krylov subspace 
	\begin{align*}
		\mathbb{G}^{j}\left({P}_{1},\ {P}_{2},\ {Q}\right) = \text{span} \{{Q},\ \mathfrak{S}_1\left({P}_1,\ {P}_2\right) {Q}, 
		\  \mathfrak{S}_2\left({P}_1,\ {P}_2\right)  {Q}, \ \ldots, \  \mathfrak{S}_j\left({P}_1,\ {P}_2\right)  {Q}\},
	\end{align*}
	where $\mathfrak{S}_j \left( {P}_1, \  {P}_2\right) =  {P}_1  \mathfrak{S}_{j-1}\left({P}_1, \ {P}_2\right) +  {P}_2  \mathfrak{S}_{j-2}\left( {P}_1,\ {P}_2\right) \ \text{for} \ j > 2, \text{with} \ \mathfrak{S}_1 \left( {P}_1, \  {P}_2\right)  = P_1 \ \text{and} \ \mathfrak{S}_2 \left( {P}_1, \  {P}_2\right) = 
	P_1^2+P_2.$
	
	For the special case of proportionally damped second-order linear systems, it has been observed in \cite{BeattieG2005} that with $\mathscr{A}=\left(s_{0}^{2}M+s_{0}D+K\right)$
	\begin{align*}
		\mathbb{G}^{j} \left( {P}_{1},\  {P}_{2},\ {Q} \right) 
		& =   \mathbb{G}^{j}\left(-\mathscr{A}^{-1} \left(2s_{0}M+D\right),\ -\mathscr{A}^{-1}M,\ \mathscr{A}^{-1}B\right), \\ 
		& = \mathbb{G}^{j} \left(-\mathscr{A}^{-1} \left( \left(2s_{0} +\alpha \right)M+\beta K \right),\ -\mathscr{A}^{-1}M,\ \mathscr{A}^{-1}B \right),\\
		&= \mathbb{K}^{j} \left( -\mathscr{A}^{-1}M,\ {\mathscr{A}^{-1}B}\right)  = \mathbb{K}^{j}\left( {P}_2,\ {Q}\right),
	\end{align*} 
	where $\mathbb{K}^{j}\left( {P}_2,\  {Q}\right)$ is the standard block Krylov subspace
	\begin{align*}
		\mathbb{K}^{j}\left( {P}_2, {Q}\right)= \text{span} \{ {Q}, \ {P}_2 {Q}, \  {P}_2^{2}  {Q}, \ \ldots, \  {P}_2^{j-1} {Q} \}.
	\end{align*}
	The reduced order system (\ref{eq:red_sys}), which matches the first $\lceil r/m \rceil$  moments of the original system 
	(\ref{eq:linear-non-para}) can be obtained by projecting (\ref{eq:linear-non-para}) with $\Pi = VV^T$ with an orthonormal matrix $V \in \mathbb{R}^{n \times r}$  whose columns span $\mathbb{K}^{j}( {P}_2, {Q}).$ 
	
	Standard efficient methods to compute the desired orthogonal basis of $\mathbb{K}^{j}\left( {P}_2,\ {Q}\right)$ are,
	e.g., the block or the global Arnoldi algorithm~\cite{Saad2003,JbiMS99,Sad93}.
	{The AIRGA algorithm}
	generates $V$ by a global Arnoldi method. Its relevant parts are given in Algorithm~\ref{Algo:AIRGA}.
	Unlike as discussed above, the AIRGA algorithm
	uses not just one expansion point, but a set of $\ell$ expansion points. This ensures a better reduced system in the entire frequency domain of interest. The method is adaptive, i.e. it automatically chooses the number of moments to be matched at each expansion point $s_i.$ This is controlled by the inner \texttt{while} loop starting at line \ref{AIRGA:while-inner}.
	The variable $j$ stores the total number of moments matched. The upper bound on {max value of $j$ or $J$ is} $\lceil r_\text{max}/m \rceil$, where $r_\text{max}$ is the maximum dimension to which we want to reduce the state variable (input from the user), and $m$ is the dimension of the input. 
	For a thorough discussion on how to determine convergence, to choose the expansion points in the inner loop as
	well as a new set of expansion points in the outer loop, see \cite{Bonin20161}.
	
	{Next, we discuss the stability analysis of using inexact linear solves in AIRGA.}
	
	{\subsection{Backward Stability Analysis}\label{sec:back_gen_con}
		Let $V$ be calculated exactly, and $f$ be the functional representation of {the exact MOR algorithm (that uses $V$ during reduction process)}. Similarly, let $\widetilde{V}$ be calculated inexactly (i.e., by a Krylov subspace solver), and $\widetilde{f}$ be the functional representation of {the inexact MOR algorithm (that uses $\widetilde{V}$ during reduction process)}. Then, from the backward stability definition, a MOR algorithm is backward stable with respect to 
		inexact linear solves if \cite{trefethen1997numerical}
		\begin{align}
			& \widetilde{f}(x) = f(\widetilde{x})   \quad \textnormal{for some $\widetilde{x}$  with } \label{eq:first_cond} \\
			& \frac{\|x- \widetilde{x}\|_{H_2 \ or \ H_{\infty}}} {\|x\|_{H_2 \ or \ H_{\infty}}} = \mathcal{O}(\|Z\|), \label{eq:secd_cond} 
		\end{align}
		where $\widetilde{x}$ is the perturbed full model corresponding to the error in the linear solves for $\widetilde{V}$ in the inexact MOR algorithm. This perturbation is denoted by $Z$. Further, $H_2$ and $H_{\infty}$ denote the standard functional norms.} 
	
	
	{
		{Here,} the function $f$ maps 
		$H(s)$ to $\hat{H}(s)$ or  
		$f(H(s))=\hat{H}(s).$
		This is represented by AIRGA 
		when a direct solver for solving the linear systems at lines \ref{algo:line-1st} and \ref{algo:line-jth} is employed {(see Algorithm \ref{Algo:AIRGA})}. This is called the exact AIRGA algorithm.}
	
	{The function $\tilde f$ maps the transfer function $H(s)$ of the original system to the transfer function of the reduced system 
		employing an iterative solver in order to solve the linear systems at lines \ref{algo:line-1st} and \ref{algo:line-jth} {in} AIRGA~(instead of a direct solver; {again see Algorithm \ref{Algo:AIRGA}}). This is denoted by $\tilde{f}(H(s))=\tilde{\hat{H}}(s)$ and is called the inexact AIRGA algorithm.} 
	
	{For our discussion, we are only interested in one outer iteration step. The matrix $V = \left[V_1, V_2, \ldots, V_J\right]$ is generated and the reduced system is computed with $V$ as in (\ref{eq:non_para_red}) (lines \ref{algo:line-V-final}-\ref{AIRGA:red_sys}). This immediately gives $f(H(s))$ and $\tilde{f}(H(s))$ when using of direct solver and iterative solver, respectively. Further, we need to assume that the choice of the expansion points is the same no matter whether iterative solves or a direct solve is used.}
	
	{Next, we analyze \eqref{eq:first_cond} and \eqref{eq:secd_cond} separately in the below two subsections.}  
	
	\subsubsection{First Condition of Stability}\label{sec:cond_stab}
	\begin{algorithm}[]
		\small
		\caption{Adaptive Iterative Rational Global Arnoldi Algorithm \cite{Bonin20161}}
		\label{Algo:AIRGA}
		\begin{algorithmic}[1]		
			\STATE Input:  \{$M,\ D,\ K,\ B,\ \mathscr{C}$,\ $r_\text{max}$; \ initial set of expansion points $S= \{s_1, \ldots, s_\ell\}$\}
			\WHILE{no convergence} \label{AIRGA:while-outer}
			\FOR  {$ \text{each}\ s_{i} \in S$} \label{AIRGA:for-1}
			\STATE {$X^{(-1)}(s_{i}) = 0$, $ h_{\pi}^{(-1)} = 0$}
			\STATE $X^{(0)}(s_{i})= \left(s_{i}^{2}M+s_{i}D+K\right)^{-1}B$, {$h_{\pi}^{(0)} = 1$} \label{algo:line-1st}
			\STATE  {Also, get a good basis of $X^{(0)}(s_i)$ via a QR decomposition}	\label{algo:qr1}		 
			\ENDFOR
			\STATE j = 1 
			\WHILE{no convergence and $j < \lceil r_\text{max} / m \rceil $} \label{AIRGA:while-inner}
			\STATE Choose an expansion point $\sigma_j \in S$; 
			{\color{black}{$\sigma_j$ = $\text{argmax}_{s_i}\|h_{\pi}^{(j-1)}\mathscr{C}^T X^{(j-1)}(s_i)\|_f$}} \label{algo:line-exp}
			\STATE {$V_{j}=X^{(j-1)}(\sigma_{j})/\|X^{(j-1)}(\sigma_{j})\|_{f}$} \label{algo:line-V}
			\FOR {$i=1,\ \ldots,\ \ell$} 
			\IF {($s_{i}==\sigma_{j}$)} 
			\STATE $X^{(j)}(s_{i})=-\left(s_{i}^{2}M+s_{i}D+K\right)^{-1}MV_{j}$, {$h_{\pi}^{(j)} = h_{\pi}^{(j-1)} \|X^{(j-1)}(s_{i})\|_{f}$}\label{algo:line-jth}
			\ELSE  
			\STATE $ \ X^{(j)}(s_{i})=X^{(j-1)}(s_{i})$,  {$h_{\pi}^{(j)} = h_{\pi}^{(j-1)}$}
			\ENDIF
			\FOR {$t=1,\ 2,\ \ldots,\ j$}\label{algo:Arnoldi-s}
			\STATE $\gamma_{t,j}(s_{i})=\text{trace}(V_{t}^{H}\cdot X^{(j)}(s_{i}))$  $X^{(j)}(s_{i})=X^{(j)}(s_{i})-\gamma_{t,j}(s_{i})V_{t}$ \label{algo:Arnoldi-a} 
			\ENDFOR \label{algo:Arnoldi-e}
			\ENDFOR
			\STATE j = j+1
			\ENDWHILE 
			\STATE Set $J = j$ and pick $\sigma_{J} \in S$
			\STATE $V_{J}=X^{(J-1)}(\sigma_{J})/||X^{(J-1)}(\sigma_{J})||_{f}$ and $V =[V_{1},\ V_{2},\ \ldots,\ V_{J}]$ \label{algo:line-V-final}
			\STATE   {Also, get a good basis of $V$ via a QR decomposition} \label{algo:qr2}
			\STATE {Compute the reduced order system matrices $\hat{M}$, $\hat{D}$ and $\hat{K}$ with $V$ as in (\ref{eq:non_para_red})} \label{AIRGA:red_sys}
			\STATE  Choose new set of expansion points $S= \{s_1, \ldots, s_\ell\}$ {using eigenvalues of the reduced system }
			\ENDWHILE
			\\Compute the reduced order system {matrices $\hat{B}$, and $\hat{\mathscr{C}}$}  with $V$ as in (\ref{eq:non_para_red}) 
		\end{algorithmic}
	\end{algorithm}	
	Consider the  linear systems for $X^{(0)}(s_{i}) \in \mathbb{R}^{n \times m}$ at line \ref{algo:line-1st} 
	\begin{align*}
		\left(s_{i}^{2}M+s_{i}D+K \right) X^{(0)} \left(s_{i} \right) = B,
	\end{align*}
	where $s_{i} \in S=\{s_{1},\ s_{2},\ \ldots,\ s_{\ell}\}$. We denote the inexactly computed solution for 
	$X^{(0)}(s_{i})$ by $\tilde{X}^{(0)}(s_{i}).$
	{Let the associated residual  be $\eta_{0i} \in \mathbb{R}^{n \times m}$ for $i=1,\ \ldots,\ \ell$. Then, the above equation is equivalent to  
		\begin{align}\label{inexact1}
			\left(s_{i}^{2}M+s_{i}D+K \right) \tilde{X}^{(0)} \left(s_{i}\right)&= B+\eta_{0i}.
		\end{align}
		{All $\tilde{X}^{(0)}(s_i)$ are used {at line \ref{algo:line-exp} for} picking the best expansion point for this {first} step, {which is denoted by} 
			$\sigma_1$ {with} $\eta_{(0)}$ {has} the corresponding residual.}
		Next, in Algorithm~\ref{Algo:AIRGA} at line  \ref{algo:line-V}, at the first iteration of the \texttt{while} loop (i.e. j=1), $\tilde{V}_{1}$ is computed as (as above, here \ $\tilde{}$  \ is added because of the inexactness)
		\begin{align}\label{v_oth}
			{\tilde{V}_{1} =  \tilde{X}^{(0)}\left(\sigma_1\right)/ \|\tilde{X}^{(0)}\left(\sigma_1\right)\|_f.}
		\end{align}   

		Further, at line \ref{algo:line-jth} in Algorithm~\ref{Algo:AIRGA}  the inexact solve gives
		\begin{align}\label{inexcat2a}
			{\left(\sigma_1^{2}M+ \sigma_1D+K \right) \tilde{X}^{(1)}\left(\sigma_1\right)= -M\tilde{V}_{1}+\eta_{1}}.
		\end{align}
		{$\tilde{X}^{(1)}(s_i)$ will be equal to $\tilde{X}^{(0)}(s_i)$,  $\forall s_i \in S = \{s_1, s_2, \ldots, s_{\ell}\} \backslash \{\sigma_1\}$}. {As above, all $\tilde{X}^{(1)}(s_i)$ are used at line \ref{algo:line-exp} for picking the best expansion point at {this} second step, which  is denoted by  $\sigma_2$ {with} $\eta_{(1)}$ {as} the corresponding residual}.	
		Next, in Algorithm~\ref{Algo:AIRGA} at line \ref{algo:line-V} after one iteration of the \texttt{while loop} (i.e. j=2), $\tilde{V}_2$ is computed as
		\begin{align}\label{eq:orth_2}
			{\tilde{V}_2 = \tilde{X}^{(1)}\left(\sigma_2\right)/\|\tilde{X}^{(1)}\left(\sigma_2\right)\|_f.}
		\end{align}
		
		Further, at line \ref{algo:line-jth} the inexact solve yields for $j = 2, \ldots, J - 1$
		\begin{align}\label{inexact2}
			{\left( \sigma_j^{2}M + \sigma_jD+K \right) \tilde{X}^{(j)}\left(\sigma_j\right)= -M\tilde{V}_{j}+\eta_{j}.}
		\end{align}
		{$\tilde{X}^{(j)}(s_i)$ will be equal to $\tilde{X}^{(j-1)}(s_i)$,  $\forall s_i \in S = \{s_1, s_2, \ldots, s_{\ell}\} \backslash \{\sigma_j\}$}. {As done for first and second step, all $\tilde{X}^{(j-1)}(s_i)$ are used  {at line \ref{algo:line-exp}} for picking the best expansion point at the $j^{th}$ step, which is denoted by  $\sigma_j$ {with} $\eta_{(J-1)}$ {as} the corresponding residual}.
		Thus, in Algorithm~\ref{Algo:AIRGA} at line \ref{algo:line-V} for $j = 3, \ \ldots, \ J-1$ and at line \ref{algo:line-V-final} for $j=J$, $\tilde{V}_{j}$ is computed as
		\begin{align}\label{eq:orth_j}
			{\tilde{V}_{j} = \tilde{X}^{(j-1)}\left(\sigma_j\right)/ \|\tilde{X}^{(j-1)}\left(\sigma_j\right)\|_f.	}
		\end{align}
		
		Finally,  $\tilde{V} = \left[\tilde{V}_1, \ \tilde{V}_2, \ \ldots, \ \tilde{V}_{J}\right]$ is set up and used to generate the reduced system (obtained by the inexact AIRGA algorithm),
		\begin{align}\label{inexact}
			\begin{split}
				& \tilde{\hat{M}} = \tilde{V}^{T}M \tilde{V},\ \tilde{\hat{D}} = \tilde{V}^{T}D \tilde{V}, \ 
				\tilde{\hat{K}} = \tilde{V}^{T} K\tilde{V}, \\ 
				& \tilde{\hat{B}} = \tilde{V}^{T}B, \ \text{and} \ {\color{black}\tilde{\hat{\mathscr{C}}}^T = \mathscr{C}^T\tilde{V}}. 
			\end{split}
		\end{align}
		This reduced order system is equivalent to $\tilde{f}\left(H(s)\right)$. 
		
		Now we have to find a perturbed original system $\tilde H\left(s\right)$, such that the exact AIRGA on it or  $f(\tilde{H}(s))$,
		will give the reduced system as obtained \big(by applying inexact AIRGA on the original full system or $\tilde{f}\left(H(s)\right)$\big). That is, find $\tilde{H}(s)$ such that $\tilde{f}\left(H(s)\right)=f(\tilde{H}(s))$. This will satisfy the first stability condition \eqref{eq:first_cond}.

		Among the many systems $\tilde{H}\left(s\right)$ one can consider here, we concentrate on those that have a constant perturbation
		$Z \in \mathbb{R}^{n \times n}$ in $K$ only. That is,
		\[ \tilde{K} =  K + Z, \ \tilde{M} = M, \ \tilde{D} = D,\ \tilde{B} = B, \ \text{and} \ \tilde{\mathscr C} = {\mathscr C}.\]
		{Although in (\ref{inexact1}), only one linear system's data is used in deciding $\tilde{V}_1$, which is $\tilde{X}^{(0)}(\sigma_1)$.  However, as mentioned earlier, all these linear solves $\tilde{X}^{(0)}(s_{i})$ are used in picking $\sigma_1$. Then, for $\tilde{H}$ we have that instead of (\ref{inexact1}), $\tilde{X}^{(0)}(s_{i})$ is the exact solution of}
		\begin{align}\label{exact1}
			{\left(s_i^{2}M+ s_iD+ \left(K+ Z\right)\right) \tilde{X}^{(0)}\left(s_i\right)= B,}
		\end{align}
		{for $i = 1, 2, \ldots, \ell$.}
		Similarly, it follows that the linear systems  (\ref{inexcat2a}) and (\ref{inexact2}) are solved exactly as
		\begin{equation}\label{exact2}
			{\left(\sigma_j^{2}M+\sigma_jD+ \left(K+Z\right)\right) \tilde{X}^{(j)}\left(\sigma_j\right)= -M\tilde{V}_{j},} 
		\end{equation}
		{for $j=1, \ldots, J-1$, where $\sigma_j$ is the expansion point picked at the $j^{th}$ step.}
		
		The final matrix $\tilde{V} = \left[\tilde{V}_1, \ \tilde{V}_2, \ldots, \ \tilde{V}_{J}\right]$ is exactly the same as before {since 
			\begin{itemize}
				\item[(a)] $\tilde{X}^{(0)}\left(\sigma_1\right)$ in (\ref{exact1}) is the same as that of (\ref{inexact1}) as well as $\tilde{X}^{(j)}\left(\sigma_j\right)$ in (\ref{exact2})  is the same as that in (\ref{inexcat2a}), (\ref{inexact2}), and
				\item[(b)] $\tilde{V}_j$ for $j = 1, \ldots, J$ are still given by (\ref{v_oth}), (\ref{eq:orth_2})  and~(\ref{eq:orth_j}). 
			\end{itemize}}
			Thus, the reduced order system \big(obtained by the exact AIRGA algorithm applied to the perturbed system $\tilde{H}$\big) is given by
			\begin{align}\label{exact}
				\begin{split}
					& \hat{\tilde{M}} = \tilde{V}^{T}\tilde{M} \tilde{V} = \tilde{V}^{T}M \tilde{V} = \tilde{\hat{M}},\\
					& \hat{\tilde{D}} = \tilde{V}^{T}\tilde{D} \tilde{V} = \tilde{V}^{T}D \tilde{V} = \tilde{\hat{D}}, \\
					&\hat{\tilde{K}} =\tilde{V}^{T} \tilde{K}\tilde{V} =  \tilde{V}^{T} \left(K+Z\right)\tilde{V} = \tilde{\hat{K}} + \tilde{V}^TZ\tilde{V}, \\ 
					& \hat{\tilde{B}} = \tilde{V}^{T}\tilde{B} = \tilde{V}^{T}B = \tilde{\hat{B}}, \ \text{and}  \\ 
					&{\color{black}\hat{\tilde{\mathscr C}}^T = \tilde{\mathscr C}^T\tilde{V} = {\mathscr C}^T\tilde{V} = \tilde{\hat{\mathscr C}}^T}.
				\end{split}
			\end{align}
			This reduced order system is equivalent to $f\left(\tilde H(s)\right)$. Obviously, this is already almost the same as $\tilde{f}\left(H(s)\right)$\Big(recall that our goal is to find $\tilde{H}\left(s\right)$ such that $\tilde{f}\left(H(s)\right)=f\left(\tilde{H}(s)\right)$\Big). 
			Thus, we need to find $Z$ such that $\hat{\tilde{K}}=\tilde{\hat{K}}$ or $ \tilde{V}^TZ\tilde{V} = 0$. 
			
			If we look at {the inexact solves in} (\ref{inexact1}), (\ref{inexcat2a}) and (\ref{inexact2}), and the corresponding perturbed {solves in} (\ref{exact1}) and (\ref{exact2}), {we find that} both are equivalent and a total {of} $\ell + J-1$ linear systems are solved. {Since} the dimension of $\tilde{V}$ is only $J$, {we further work with} only those linear systems {that} form our  $\tilde{V}$ {and} ignore the remaining systems. {Putting} all these linear systems together {we} get 
			{
				\begin{align}\label{eq:Orth_X}
					Z \ \mathbf{X} &= \eta,
				\end{align}
				where $\mathbf{X}$ is formed by stacking {the relevant} block columns of $X^{(j)}(\sigma_j)$ {or} $\mathbf{X}= \left[\tilde{X}^{(0)}\left(\sigma_1\right), \ \tilde{X}^{(1)}\left(\sigma_2\right), \ \ldots, \ \tilde{X}^{(J-1)}\left(\sigma_{J}\right) \right]$; similarly, after stacking {the relevant} block columns of $\eta_{j}$ together we {get} $\eta = \left[-\eta_{(0)},\ \ldots,  -\eta_{(J-1)}\right]$.}

			
			{
				In the above equation,  we can replace $\mathbf{X}$ in-terms of $\tilde{V}$ by using (\ref{v_oth}), (\ref{eq:orth_2}), and (\ref{eq:orth_j}). 
				That is, (\ref{eq:Orth_X}) can be rewritten as
				\begin{align}\label{delta_K1}
					Z \tilde{V} \mathcal{D}_X^{-1} =\eta \quad or \ \
					Z \tilde{V} = \eta \mathcal{D}_X,
				\end{align}
				where 
			} 
			$\mathcal{D}_X = \begin{bmatrix}
			\frac{1}{\left\|\tilde{X}^{(0)}\left(\sigma_1\right)\right\|_f} & 0 & 0 & \cdots & 0\\
			0 & \frac{1}{\left\|\tilde{X}^{(1)}\left(\sigma_2\right)\right\|_f} & 0 & \cdots &  0 \\
			\vdots & \vdots & \vdots & \ddots & \vdots \\ 
			& & & & \\
			0 & 0 & 0 & \cdots & \frac{1}{\left\|\tilde{X}^{(J-1)}\left(\sigma_J\right)\right\|_f}\\
			\end{bmatrix}.$} 
		
		Multiplying $\tilde{V}^T$ from the left side of (\ref{delta_K1}), we get
		\begin{align}\label{V_oth_Z}
			\tilde{V}^{T}Z \tilde{V} = \tilde{V}^{T} \eta \mathcal{D}_X.
		\end{align}
		{Assume that we are using a Ritz-Galerkin based iterative solver. Here, the solution space of the linear systems is orthogonal to the corresponding residuals, i.e. 
			{$\tilde{V}_{1} \perp {\eta}_{(0)}, \ \tilde{V}_{2} \perp {\eta}_{(1)},\ \ldots \ ,\text{and} \ \tilde{V}_{J}  \perp  {\eta}_{(J-1)}$}~\cite{van2003iterative}. Hence, 
			{\begin{align}\label{oth_matrix}
					\begin{split}
						\tilde{V}^{T} \eta &= - \begin{bmatrix}
							\tilde{V}_{1}^{T} \\
							\tilde{V}_{2}^{T} \\
							\vdots &  \\
							\tilde{V}_{J-1}^{T} \\
							\tilde{V}_{J}^{T}
						\end{bmatrix}\begin{bmatrix}
						{\eta}_{(0)} &
						{\eta}_{(1)} &
						\ldots &  
						{\eta}_{(J-1)}   
					\end{bmatrix},   \\ &= - \begin{bmatrix}  0 &  \tilde{V}^{T}_{1}{\eta}_{(1)} & \ldots &  \tilde{V}^{T}_{1}{\eta}_{(J-2)}&  \tilde{V}^{T}_{1}{\eta}_{(J-1) }  \\    \tilde{V}^{T}_{2}{\eta}_{(0)} &  0 & \ldots &  \tilde{V}^{T}_{2}{\eta}_{(J-2) } & \tilde{V}^{T}_{2}{\eta}_{(J-1) }   \\   \vdots & \vdots & \vdots & \vdots &\vdots\\
					\tilde{V}^{T}_{J-1}{\eta}_{(0)}  & \tilde{V}^{T}_{J-1}{\eta}_{(1)} & \ldots & 0 &  \tilde{V}^{T}_{J-1}{\eta}_{(J-1) }\\  \tilde{V}^{T}_{J}{\eta}_{(0)} &\tilde{V}^{T}_{J}{\eta}_{(1)} & \ldots & \tilde{V}^{T}_{J}{\eta}_{(J-2) } & 0 \end{bmatrix}.
			\end{split}
		\end{align}}}
		{ 
			Our goal here is to make the right hand side of the above equation equal to zero. The upper triangular part of the above matrix is zero if we have the following orthogonalities:
			{ \begin{align}\label{uppr_trig_oth}
					\begin{split}
						\begin{bmatrix}
							\tilde{V}_1
						\end{bmatrix} &\perp {\eta}_{(1) },  \\
						\begin{bmatrix}
							\tilde{V}_1 \ \tilde{V}_2  
						\end{bmatrix} &\perp {\eta}_{(2) },  \\
						\vdots &  \\
						\begin{bmatrix}  
							\tilde{V}_1 \ \tilde{V}_2 \ \tilde{V}_3 \ \ldots   \ \tilde{V}_{J-2}
						\end{bmatrix} &\perp {\eta}_{(J-2) },  \\ 
						\begin{bmatrix}  
							\tilde{V}_1 \ \tilde{V}_2 \ \tilde{V}_3 \ \ldots  \tilde{V}_{J-2} \ \tilde{V}_{J-1}
						\end{bmatrix} &\perp {\eta}_{(J-1) }. 
					\end{split}
				\end{align}}
				%
				Similarly, for the lower triangular part of the above matrix to be zero we need the following orthogonalities:
				{\begin{align}\label{lower_trig_oth}
						\begin{split}
							\tilde{V}_{2} \perp & 
							\begin{bmatrix}
								{\eta}_{(0)}
							\end{bmatrix},  \\ 
							\tilde{V}_{3} \perp & \begin{bmatrix}
								{\eta}_{(0)} \ {\eta}_{(1) }
							\end{bmatrix},   \\ 
							\vdots &  \\ 
							\tilde{V}_{J-1} \perp & \begin{bmatrix}  
								{\eta}_{(0)} \ {\eta}_{(1 )} \ \ldots \ {\eta}_{(J-3) }
							\end{bmatrix},  \\ 
							\tilde{V}_{J} \perp & \begin{bmatrix}  
								{\eta}_{(0)} \ {\eta}_{(1) } \ \ldots \ {\eta}_{(J-3) } \ {\eta}_{(J-2) }
							\end{bmatrix}.  
						\end{split}
					\end{align}}
					At the first glance, there seem to be two problems in achieving the above discussed orthogonalities in an iterative solver. One is the amount of code changes to be done. The other is the extra cost associated at every iterative step of the solver, which may undermine the benefit of using an iterative solver itself.
					In Section~\ref{sec:implementation_rcg}, we show that both these issues can be easily resolved by using a recycling variant of the underlying iterative solver (briefly summarized below).
					
					While solving a sequence of linear systems, if the consecutive systems do not change much, then some information can be reused from solving one linear system to solving the next. In the context of Krylov based iterative linear solvers, this information is in the form of the generated Krylov subspace. The process of  reusing Krylov subspaces from one linear system to the next is termed as ``recycling"~\cite{Michael2006, Wang2007, Ahuja2012, kapil2012}.
					
					A subset of $\tilde{V}$'s and $\eta$'s of  (\ref{uppr_trig_oth}) and (\ref{lower_trig_oth}) can be used to span a recycle space, leading to almost no code changes in the recycling variant of the underlying iterative solver. In some cases, this choice of the recycle space can actually accelerate the convergence of the next linear system in the sequence. In case when this recycle space deteriorates the convergence of the next linear system, this behaviour is bounded. In the numerical experiments section (Section~\ref{sec:Numerical_results}), we support both these conjectures (acceleration and deterioration of the convergence of iterative linear solvers) with multiple examples.

					Therefore, after applying (\ref{oth_matrix}), (\ref{uppr_trig_oth}) and (\ref{lower_trig_oth}) to~(\ref{V_oth_Z}),
					we get  $\tilde{V}^{T} Z \tilde{V} =0$. 
					Thus, $\hat{\tilde{K}} = \tilde{\hat{K}}$ or
					%
					%
					%
					\begin{align*}
						\tilde{f}\left(H\left(s\right)\right)=f\left(\tilde{H}\left(s\right)\right)=\tilde{\hat{H}}\left(s\right),
					\end{align*}  	
					$\text{where} \ H\left(s\right)= \mathscr{C}^T \left(s^{2}M+sD+K\right)^{-1}B,\\ \tilde{H}\left(s\right) = \mathscr{C}^T \left(s^{2}M+sD+  \left(K+ Z\right)\right)^{-1}B,  \text{and}  \\ 
					\tilde{\hat{H}}\left(s\right) = \hat{\mathscr{C}}^T\left(s^{2}\tilde{\hat{M}}+s\tilde{\hat{D}}+\tilde{\hat{K}}\right)^{-1}\tilde{\hat{B}} =
					\hat{\mathscr{C}}^T  \left(s^{2}\hat{\tilde{M}}+s\hat{\tilde{D}}+\hat{\tilde{K}}\right)^{-1}\hat{\tilde{B}}$.
					Thus, we satisfy the first condition of stability.
					
					\subsubsection{Second Condition for Stability}\label{sec:second-cond}
					According to the second condition of stability, given in (\ref{eq:secd_cond}), the difference between the unperturbed (original) full system and the perturbed full system should be of the order of the perturbation~\cite{trefethen1997numerical}. These errors are measured in the commonly used norms as below.
					\begin{align*}
						& {H}_{2}-\textbf{norm} \qquad \qquad \| H - G\|_{H_2} = \frac{1}{2\pi} \int_{-\infty}^\infty \| H(\imath \omega) - G(\imath \omega)\|_f d\omega, \\
						& {H}_{\infty}-\textbf{norm} \qquad \qquad  \| H - G\|_{H_\infty} = \max_{\omega \in \mathbb{R}} \| H(\imath \omega) - G(\imath \omega)\|_2,
					\end{align*}
					where the transfer functions $H$ and $G$ belong to systems with the same input and output dimension. 
					Theorem 4.3 from~\cite{Beattie20122916} helps in giving the desired result.
					\\ \quad \\
					{\it Theorem 1:}
						$ \ \text{If}  \ \  \|Z\|_2 < \frac{1}{\|A\left(s\right)^{-1}\|_{H_{\infty}}} \ \ \text{then} $
						\begin{align}\label{eq:send_cond}
							\|H(s)-\tilde{H}(s)\|_{H_2} \le \frac{\|A(s)^{-1}B\|_{H_{\infty}} \|{\mathscr C}^TA(s)^{-1}\|_{H_{2}}}{1-\|A(s)^{-1}\|_{H_{\infty}}\|Z\|_2}\|Z\|_2, 
						\end{align} 
						where $A(s) = \left(s^2 M + s D+ K\right)$, 
						{$H(s) = \mathscr{C}^T A(s)^{-1} {B}$ and  $\tilde{H}(s) = \mathscr{C}^T \big(A(s) + Z\big)^{-1} {B}$}.
					{\it Proof: We use Theorem 4.3 from \cite{Beattie20122916} for this.}
					
					If $\|Z\|_2 < 1$ and $\|A(s)^{-1}\|_{H_\infty} < 1$, then we have $\|A(s)^{-1}\|_{H_\infty}  \|Z\|_2 < 1$, and hence,
					\begin{align}\label{eq:bound_sec_cond}
						\frac{1}{1- \|A(s)^{-1}\|_{H_\infty} \|Z\|_2} < \frac{1}{1- \|A(s)^{-1}\|_{H_\infty}}.
					\end{align} 
					
					Substituting (\ref{eq:bound_sec_cond}) in (\ref{eq:send_cond}) we get
					\begin{align}\label{eq:sec_con_bound}
						\begin{split}
							\frac{\| H(s)-\tilde{H} ( s  ) \|_{H_{2}} }{\|H(s)\|_{H_2}} 
							\le & \ \frac{\|A(s)^{-1}B\|_{H_{\infty}} \|{\mathscr C}^TA(s)^{-1}\|_{H_{2}}}{\|H(s)\|_{H_2}} \cdot \frac{1}{1- \|A(s)^{-1}\|_{H_\infty}} \cdot \|Z\|_2 \\
							= & \ \mathcal{O}(\| Z \|_2).
						\end{split}
					\end{align}
					{This satisfies the second condition of stability.}
					
					{The next theorem summarizes our complete stability analysis}
					\\ \quad \\
					{\it Theorem 2:}
						If the linear systems arising in the AIRGA algorithm are solved by 
						
						(a) Ritz-Galerkin based solver (i.e. the residual is orthogonal to the generated Krylov subspace), 
						
						(b)	 the extra orthogonalities given by (\ref{uppr_trig_oth}) and~(\ref{lower_trig_oth}) are satisfied by such a solver,
						
						(c) $A\left(s\right)$ as defined in Theorem 1 is invertible and $\|A(s)^{-1}\|_{H_\infty} < 1 $,
						
						(d)	$Z$ is given by (\ref{eq:Orth_X}) exists and  $\|Z\|_2 < 1$, 
						
						then the AIRGA algorithm is backward stable with respect to the inexact linear solves. 
					If we look at Algorithm \ref{Algo:AIRGA}, besides the linear solves at lines \ref{algo:line-1st} and \ref{algo:line-jth}, we are also concerned about the construction of $V$ from $X$ (since $V$ gives us our reduced system). There are three places in code  where $X$ is modified further to obtain $V$ rather than just normalizing $X$ to $V$ ({on lines \ref{algo:line-V} and \ref{algo:line-V-final}}). First is at line \ref{algo:qr1}, where a QR decomposition of $X$ is done. Second, an Arnoldi iteration on $X$ is done at lines \ref{algo:Arnoldi-s}--\ref{algo:Arnoldi-e}. Finally, the QR decomposition of $V$ is done at line \ref{algo:qr2}. All these code changes are nothing but an attempt to get a good basis of $X$ and $V,$ {which have negligible effect on our analysis. Hence, for ease of exposition, we ignore them.} 
					
					{Next, we discuss how to satisfy the backward stability conditions given by Theorem 2}
					\section{Satisfying Backward Stability Conditions} \label{sec:Sat_Back_stab_cond}
					{In this section,  we analyze the hypothesis of Theorem 2 as so to achieve a backward stable AIRGA, we mostly focus on conditions (a) and (b) and not (c) and (d). Condition (c) cannot be worked upon much because it is dependent on the expansion points and the input dynamical system. Condition (d) does not create much challenges because, as discussed in the next section, perturbation is directly proportional to the residuals, which can be iteratively controlled.}
					
					{From condition (a) of Theorem 2 we know that} 
					we need to use a Ritz-Galerkin based method for solving the underlying linear systems {in AIRGA}. The Conjugate Gradient (CG) method is one of the most popular solver of such a type. The CG method is mainly used for solving Symmetric Positive Definite (SPD) linear systems. For solving non-symmetric linear systems, Full Orthogonalization Method (FOM)~\cite{JbiMS99,YOUNG1980} is the one that is based upon the Ritz-Galerkin theory.
					
					In this work, we focus on the CG method, {and hence,} 
					in the results section, we take models that lead to SPD linear systems in the AIRGA algorithm. 
					FOM method {can be similarly used}.
					
					Next, we first discuss how to change the theory of the CG method such that {condition (b) of Theorem 2 or} the extra orthogonalities, (\ref{uppr_trig_oth})-(\ref{lower_trig_oth}), are satisfied (in Section~\ref{sec:Ach_ext_orth}). Further, we describe how the recommended changes can be easily implemented (in Section~\ref{sec:implementation_rcg}).  
					
					{
						\subsection{Achieving Extra Orthogonalities}\label{sec:Ach_ext_orth}
						The CG method consists of two components. One is the Lanczos algorithm that gives a good basis of the generated Krylov subspace. The other is the Ritz-Galerkin projection to obtain solution estimates from this subspace.
						The orthogonalities in (\ref{lower_trig_oth}) can be achieved by modifying the Lanczos algorithm (discussed in Section~\ref{sec:Lanczos_exta_orth}), and those in (\ref{uppr_trig_oth}) can be achieved by changing the Ritz-Galerkin projection (discussed in Section~\ref{sec:Ritz-galerkin_extra_orth}). 
						
						\subsubsection{Adapting the Lanczos Process}\label{sec:Lanczos_exta_orth}
						Assume we are trying to solve the  linear system of the form
						\begin{align}\label{eq:gen_eq}
							Ax = b,
						\end{align}
						where $A \in  \mathbb{C}^{n \times n}$ and $b  \in \mathbb{C}^{n}$. Let $x_0$ be the initial solution vector with $r_0 = b - A x_0$ as the corresponding residual. 
						The Lanczos algorithm computes a good basis of the generated Krylov subspace involving $A$ and $r_0$ as~\cite{Saad2003}    
						\begin{align}\label{eq:gen_eq_Lan_basis}
							\begin{split}
								w_{k+1} \ \in \ \mathbb{K}^k\left(A,r_0\right) \equiv \text{span}\{r_0, Ar_0, A^2r_0, \cdots, A^{k-1}r_0\}
								\\
								s.t. \qquad  \ w_{k+1} \perp \begin{bmatrix}
									w_1 & w_2 & \cdots & w_k
								\end{bmatrix},
							\end{split}
						\end{align}
						where $w_{k+1}$ is the Lanczos vector at the $\left(k+1 \right)^{th}$ iterative step and $w_1 = r_0/\|r_0\|$\footnote{Here, the first equation of (\ref{eq:gen_eq_Lan_basis}) is implemented using 
							\begin{align*}
								w_{k+1} = Aw_k-c_1w_1-c_2w_2- \ldots - c_{k-1}w_{k-1}-c_kw_k.
							\end{align*}
							Finally, the second equation of 	(\ref{eq:gen_eq_Lan_basis}) gives us $c_1,c_2, \ldots, c_k$. For a complete derivation of this, please see chapter $5$ of \cite{van2003iterative}.}. 
						Now, assume we are carrying some residual vector $\tilde{r}$ from another linear system, which we need to make orthogonal to the final solution of (\ref{eq:gen_eq}). Then, the Lanczos algorithm above would consist of the following procedure:
						\begin{align*}
							w_{k+1} \in \mathbb{K}^{k}\left(A,r_0\right) \\
							s.t. \qquad w_{k+1} \perp \left[
							w_1 \ w_2 \ \cdots w_k \ \underbrace{\tilde{r}}
							\right].
						\end{align*}
						
						Recall from the previous sections that in AIRGA, the first set of linear systems to be solved iteratively are given by \eqref{inexact1}. As mentioned in the paragraph between \eqref{inexact1}--\eqref{v_oth}, the expansion point chosen is $\sigma_1$, and hence, in the linear system playing a role in our stability analysis is
						\begin{align}\label{eq:inexact1}
							\left(\sigma_1^{2}M+\sigma_1D+K \right) \tilde{X}^{(0)}\left(\sigma_1\right) = B + \eta_{(0)}.
						\end{align}
						
						Next, we need to iteratively solve (\ref{inexcat2a}), i.e.
						\begin{align}\label{eq:inexcat2a}
							\left(\sigma_1^{2}M+\sigma_1 D+K\right) \tilde{X}^{(1)}\left(\sigma_1\right)= -M\tilde{V}_{1} +\eta_{1}.
						\end{align}
						Here, we need a good basis of the Krylov subspace involving {the coefficient matrix} {$\mathcal{K}_{1} = \left(\sigma_1^2 M + \sigma_1 D + K \right)$} 
						and {$\left(\eta_1\right)_0$}
						, which is the initial residual of (\ref{eq:inexcat2a}). Hence, the Lanczos algorithm here would consist of the following procedure:
						\begin{align}\label{eq:lanczos_second_AIRGA}
							\begin{split}
								\left(w_1\right)_{k+1} \in \mathbb{K}^k \left(\mathcal{K}_1, \left(\eta_1\right)_0 \right) \\
								s.t. \qquad \left(w_1\right)_{k+1} \perp \left[
								\left(w_1\right)_1 \ \left(w_1\right)_2 \ \cdots \ \left(w_1\right)_{k}
								\right],
							\end{split}
						\end{align}	
						where  {$\left(w_1\right)_{k+1}$} 
						is the Lanczos vector at the $(k+1)^{th}$ iterative step and {$\left(w_1\right)_1 = \left(\eta_1\right)_0/\|\left(\eta_1\right)_0\|$}. 
						{At this stage it is not clear if the solution of \eqref{eq:inexcat2a}, i.e. $\tilde{X}^{(1)}(\sigma_1)$, would be used to form $\tilde{V}_2$ or not (see the discussion between \eqref{inexcat2a}--\eqref{eq:orth_2}; equivalently lines \ref{algo:line-exp}--\ref{algo:line-V} of Algorithm \ref{Algo:AIRGA}). However, to avoid repeating solving \eqref{eq:inexcat2a} incase its solution is used to form $\tilde{V}_2$, we adapt the Lanczos procedure given by \eqref{eq:lanczos_second_AIRGA} as}
						
						{\begin{align}\label{eq:ext_Lan_1st}
								\begin{split}
									\left(w_{1}\right)_{k+1} \in \mathbb{K}^k \left(\mathcal{K}_1, \left(\eta_1\right)_0 \right) \\
									s.t. \qquad \left(w_1\right)_{k+1} \perp \left[
									\left(w_1\right)_1 \ \left(w_1\right)_2 \ \cdots \ \left(w_1\right)_k \ \underbrace{\eta_{(0)}}
									\right],
								\end{split}
							\end{align}}
							where $\eta_{(0)}$ is the final residual obtained after solving (\ref{eq:inexact1}) iteratively.

							{Next, the expansion point $\sigma_2$ is chosen (as above, see paragraph between \eqref{inexcat2a}--\eqref{eq:orth_2} or lines \ref{algo:line-exp}--\ref{algo:line-V} of Algorithm \ref{Algo:AIRGA}). If $\sigma_2$ turns to be equal to $\sigma_1$ (i.e. $\sigma_2 = \sigma_1$), then $\tilde{V}_2 = \tilde{X}^{(1)} (\sigma_1)/\|\tilde{X}^{(1)} (\sigma_1)\|_f$, and we would be satisfied the {\it first} orthogonality of \eqref{lower_trig_oth}, i.e. $\tilde{V}_2 \perp [\eta_{(0)} ]$.}
							
							{If $\sigma_2$ turns to be not equal to $\sigma_1$ (say $\sigma_2 = s_i \neq \sigma_1$), then $\tilde{V}_2 = \tilde{X}^{(1)} (s_i)/\|\tilde{X}^{(1)} (s_i)\|_f = \tilde{X}^{(0)} (s_i)/\|\tilde{X}^{(0)} (s_i)\|_f$ with $\tilde{X}^{(0)}(s_i)$ given by \eqref{inexact1} or
								\begin{align}\label{eq:inexact1_rep}
									\left(s_{i}^{2}M+s_{i}D+K \right) \tilde{X}^{(0)} \left(s_{i}\right)&= B+\eta_{0i},
								\end{align}
								which we would have already solved once. Hence, to satisfy the {\it first} orthogonality of \eqref{lower_trig_oth} or $\tilde{V}_2 \perp [\eta_{(0)}]$, we would need to resolve \eqref{inexact1} or \eqref{eq:inexact1_rep} by adapting its Lanczos process as given in \eqref{eq:ext_Lan_1st}. That is, carry extra $\eta_{(0)}$ in its Krylov subspace.
							}
							
							Next, we need to iteratively solve (\ref{inexact2}) for $j = 2$, i.e.
							\begin{align}\label{eq:inexact3a}
								\left(\sigma_2^{2}M+\sigma_2D+K\right) \tilde{X}^{(2)}\left(\sigma_2\right)= -M\tilde{V}_{2} +\eta_{2}.
							\end{align}
							Here, we need a good basis of the Krylov subspace involving {the coefficient matrix} {$\mathcal{K}_2 = \left(\sigma_2^{2}M+ \sigma_2D+K \right) $ and $\left(\eta_2\right)_0$}
							, which is the initial residual of (\ref{eq:inexact3a}). Hence, the Lanczos algorithm here would consist of the following procedure:  
							\begin{align}\label{eq:lanczos_second_Lystem}
								\begin{split}
									\left(w_2\right)_{k+1}  \in \mathbb{K}^k \left(\mathcal{K}_2,\left(\eta_2\right)_0 \right)  \\
									s.t. \qquad \left(w_2\right)_{k+1}  \perp \left[
									\left(w_2\right)_1 \ \left(w_2\right)_2 \ \cdots \ \left(w_2\right)_k 
									\right],
								\end{split}
							\end{align}
							where  {$\left(w_2\right)_{k+1}$} is the Lanczos vector at the $(k+1)^{th}$ iterative step and {$\left(w_2\right)_1 = \left(\eta_2\right)_0/\|\left(\eta_2\right)_0\|$}. 
							
							{As earlier, at this stage it is not clear if the solution of \eqref{eq:inexact3a}, i.e. $\tilde{X}^{(2)} (\sigma_2)$, would be used to form $\tilde{V}_3$ or not (see the discussion between \eqref{inexact2}--\eqref{eq:orth_j}; equivalently lines \ref{algo:line-exp}--\ref{algo:line-V} of Algorithm \ref{Algo:AIRGA}). However to avoid repeating solving \eqref{eq:inexact3a} incase its solution is used to form $\tilde{V}_3$, we adapt the Lanczos procedure given by \eqref{eq:lanczos_second_Lystem} as
							}
							{\begin{align}\label{eq:ext_orth_lanc_2nd}
									\begin{split}
										\left(w_2\right)_{k+1}  \in \mathbb{K}^k \left(\mathcal{K}_2,\left(\eta_2\right)_0\right)\\
										s.t. \qquad \left(w_2\right)_{k+1}  \perp \left[
										\left(w_2\right)_1 \ \left(w_2\right)_2 \ \cdots \ \left(w_2\right)_k \ \underbrace{\eta_{(0)} \ \eta_{(1)}}
										\right],
									\end{split}
								\end{align}}
								where $\eta_{(1)}$ is the final residual obtained after solving (\ref{eq:inexcat2a}) iteratively.
								
								{Next, the expansion point $\sigma_3$ is chosen (as above, see paragraph between \eqref{inexact2}--\eqref{eq:orth_j} or equivalently lines \ref{algo:line-exp}--\ref{algo:line-V} of Algorithm \ref{Algo:AIRGA}). If $\sigma_3$ turns to be equal to $\sigma_2$ (i.e. $\sigma_3 = \sigma_2$), then $\tilde{V}_3 = \tilde{X}^{(2)} (\sigma_2)/\|\tilde{X}^{(2)} (\sigma_2)\|_f$, and we would have satisfied the {\it second} set of orthogonalities of \eqref{lower_trig_oth}, i.e. $\tilde{V}_3 = [\eta_{(0)} \ \eta_{(1)}]$.}
								
								{If $\sigma_3$ turns to be not equal to $\sigma_2$ (i.e. $\sigma_3 \neq \sigma_2$), then it may be equal to $\sigma_1$ or some other $s_i$. Incase $\sigma_3 = \sigma_1 \neq \sigma_2$, than
									$\tilde{V}_3 = \tilde{X}^{(2)} (\sigma_1)/\|\tilde{X}^{(2)} (\sigma_1)\|_f = \tilde{X}^{(1)} (\sigma_1)/\|\tilde{X}^{(1)} (\sigma_1)\|_f$ with $\tilde{X}^{(1)}(\sigma_1)$ given by \eqref{inexcat2a} or \eqref{eq:inexcat2a}, which we would have already solved once. Hence, to satisfy the {\it second} set of orthogonalities of \eqref{lower_trig_oth} or $\tilde{V}_3 \perp [\eta_{(0)} \ \eta_{(1)}]$, we would need to resolve \eqref{inexcat2a} or \eqref{eq:inexcat2a} by adapting its Lanczos process as given in \eqref{eq:ext_orth_lanc_2nd}. That is, carry extra $\eta_{(0)}$ and  $\eta_{(1)}$ in its Krylov subspace.
								}

								{Incase $\sigma_3 = s_i$ with $s_i  \neq \sigma_1$, and $s_i \neq \sigma_2$, than
									$\tilde{V}_3 = \tilde{X}^{(2)} (s_i)/\|\tilde{X}^{(2)} (s_i)\|_f = \tilde{X}^{(1)} (s_i)/\|\tilde{X}^{(1)} (s_i)\|_f = \tilde{X}^{(0)} (s_i)/\|\tilde{X}^{(0)} (s_i)\|_f$ with $\tilde{X}^{(0)}(s_i)$ 
									again given by \eqref{inexact1} or}
								\begin{align}\label{inexact1-resol}
									\left(s_{i}^{2}M+s_{i}D+K \right) \tilde{X}^{(0)} \left(s_{i}\right)&= B+\eta_{0i},
								\end{align} 
								{which we would have already solved once. Hence, to satisfy the {\it second} set of orthogonalities of \eqref{lower_trig_oth} or $\tilde{V}_3 \perp [\eta_{(0)} \ \eta_{(1)}]$, we would need to resolve \eqref{inexact1} or \eqref{inexact1-resol} by adapting its Lanczos process as given in \eqref{eq:ext_orth_lanc_2nd}. That is, carry extra $\eta_{(0)}$ and  $\eta_{(1)}$ in its Krylov subspace.}

								We need to repeat a similar procedure for (\ref{inexact2}) for all $j = 3, \ldots, J-1$. 
								\subsubsection{Adapting the Ritz-Galerkin Projection}\label{sec:Ritz-galerkin_extra_orth}
								
								Recall that if we were trying to solve the linear system given in~(\ref{eq:gen_eq}) by the CG method, then~(\ref{eq:gen_eq_Lan_basis}) gives a good basis of the generated Krylov subspace. The solution update {at the $k^{th}$ iterative step is given as \cite{Saad2003}}
								\begin{align}\label{eq:gen_eq_solution}
									x_k = x_0 + \zeta_k,
								\end{align}
								where $\zeta_k  = W_k  y_k$ and $W_k = \left[w_1 \ w_2 \ \ldots \ w_k\right]$ with the columns of this matrix given by (\ref{eq:gen_eq_Lan_basis}). In the CG method, this $y_k$ is defined by a Ritz-Galerkin projection
								\begin{align}
									r_k  \perp W_k,
								\end{align}
								where $r_k = b - A(x_0+\zeta_k) = b - A(x_0 + W_k y_k) = r_0 - A  W_k   y_k$. 		
								Now, assume we are carrying some solution vector $\tilde{x}$ from another linear system, which we need to make orthogonal to the final residual of~(\ref{eq:gen_eq}). Then, the Ritz-Galerkin projection as above would  consists of the following procedure:
								\begin{align}
									r_k  \perp \begin{bmatrix}
										W_k & \underbrace{\tilde{x}}
									\end{bmatrix}.
								\end{align}

								Let us now look at the second linear system to solve in the AIRGA algorithm, i.e. (\ref{eq:inexcat2a}). For this, a good basis of the generated Krylov subspace is given by (\ref{eq:lanczos_second_AIRGA}). To find the solution vector here, the Ritz-Galerkin projection is defined as 
								\begin{align}\label{eq:second_resid_AIRGA}
										\left(\eta_1\right)_k \perp \left(W_1\right)_k,
									\end{align}
									\sloppy where {$\left(\eta_1\right)_k$} 
									is the residual of (\ref{eq:inexcat2a}) at the $k^{th}$ iterative step and 
									{$\left(W_1\right)_k = \big[\left(w_1\right)_1 \ \left(w_1\right)_2 \ \cdots \ \left(w_1\right)_k \big]$}
									with the columns of this matrix given by (\ref{eq:lanczos_second_AIRGA}). Note that, {as earlier}, {$\eta_1$} 
									is the final residual of (\ref{eq:inexcat2a}) (at convergence of CG). 
									
									{As earlier, at this stage it is not clear if the residual of \eqref{eq:inexcat2a}, i.e. $\eta_1$ would be the residual we care, about, i.e. $\eta_{(1)}$. These two residuals map to the fact whether $\tilde{X}^{(1)}(\sigma_1)$ would be used to form $\tilde{V}_2$ or not (again see the discussion between \eqref{inexcat2a}--\eqref{eq:orth_2}; equivalently lines \ref{algo:line-exp}--\ref{algo:line-V} of Algorithm \ref{Algo:AIRGA}). However, to avoid repeating solving \eqref{eq:inexcat2a} incase its solution is used to form $\tilde{V}_2$, we adapt the projection given by \eqref{eq:second_resid_AIRGA} as 
									}
									{\begin{align}\label{eq:second_resid_AIRGA1}
											\left(\eta_1\right)_k \perp \left[
											\left(W_1\right)_k \ \underbrace{\tilde{V}_1}
											\right],
										\end{align}}	
										{where $\tilde{V}_1$ is given by (\ref{v_oth}).} 
										
										{Next, the expansion point $\sigma_2$ is chosen (as above, see paragraph between \eqref{inexcat2a}--\eqref{eq:orth_2}; equivalently lines \ref{algo:line-exp}--\ref{algo:line-V} of Algorithm \ref{Algo:AIRGA}). If $\sigma_2$ turns to be equal to $\sigma_1$ (i.e. $\sigma_2 = \sigma_1$), then $\tilde{V}_2 = \tilde{X}^{(1)} (\sigma_1)/\|\tilde{X}^{(1)} (\sigma_1)\|_f$, and we would be satisfied the {\it first} orthogonality of \eqref{uppr_trig_oth}, i.e. $\tilde{V}_1 \perp [\eta_{(1)} ]$.}
										
										{If $\sigma_2$ turns to be not equal to $\sigma_1$ (say $\sigma_2 = s_i \neq \sigma_1$), then $\tilde{V}_2 = \tilde{X}^{(1)} (s_i)/\|\tilde{X}^{(1)} (s_i)\|_f = \tilde{X}^{(0)} (s_i)/\|\tilde{X}^{(0)} (s_i)\|_f$ with $\tilde{X}^{(0)}(s_i)$ given by \eqref{inexact1} or
											\begin{align}\label{eq:inexact1_rep1}
												\left(s_{i}^{2}M+s_{i}D+K \right) \tilde{X}^{(0)} \left(s_{i}\right)&= B+\eta_{0i},
											\end{align}
											which we would have already solved. Hence, to satisfy the {\it first} orthogonality of \eqref{uppr_trig_oth} or $\tilde{V}_1 \perp [\eta_{(1)}]$, we would need to resolve \eqref{eq:inexact1_rep1} by adapting its projection as given in \eqref{eq:second_resid_AIRGA1}. That is, carrying extra $\tilde{V}_1$ in its Krylov subspace.
										}
										
										{Similarly, all the other orthogonalities of  \eqref{uppr_trig_oth} can be achieved. As mentioned earlier, the use of recycling variant of CG helps us avoid the cumbersome code changes, and this discussed next.}

										\subsection{Implementation}\label{sec:implementation_rcg}

										{
											Developing the CG algorithm that is based upon the adapted Lanczos process and the adapted Ritz-Galerkin projection is doable. However, developing its efficient implementation involving  standard two/ three term recurrences is non-trivial. Also, as the sequence number of the linear system increase (i.e. $j$ gets larger), the number of orthogonalizations to be done also increase linearly.

											As briefly discussed in Section \ref{sec:cond_stab}, using a recycling CG (RCG)~\cite{Michael2005, Parks2010} helps alleviate both these problems. Hence, in this subsection we \textit{first} discuss the idea behind RCG. \textit{Second} we describe how to use RCG so as to easily achieve the earlier described extra orthogonalities. We do this with no code changes to the existing algorithm. Here, we also discuss the extra computational cost of such an implementation.   
											
											Assume that we want to solve the linear system in (\ref{eq:gen_eq}). Also assume that the recycle space is in the form of span \{U\}, where columns of $U \in \mathbb{R}^{n \times k}$ are linearly independent. If $x_{-1}$ is the initial guess for (\ref{eq:gen_eq}) and $r_{-1} = b - Ax_{-1}$ is the corresponding residual, then the projected initial guess $x_0$ is defined as~\cite{Michael2005, Parks2010}
											\begin{align*}
												x_0 = x_{-1} + U\left(U^TAU\right)^{-1}U^Tr_{-1},
											\end{align*}
											with the corresponding residual $r_0 = b -Ax_0$.
											
											At the $k^{th}$ iterative step, the Lanczos process involves \cite{Saad2000defcg}
											\begin{align*}
												w_{k+1} \in \mathbb{K}^k(A,U,r_0) \equiv \text{span}\{U, r_0, Ar_0, A^2 r_0, \cdots, A^{k-1}r_0\} \\
												s.t. \qquad w_{k+1} \perp \left[U \ w_1 \ w_2 \ \cdots w_k \right],
											\end{align*}
											where $w_{k+1}$, as earlier, is the $\left(k+1\right)^{th}$ Lanczos vector and $w_1 = r_0/\|r_0\|$. 
											The Ritz-Galerkin projection here is as follows:
											\begin{align*}
												r_k \perp \mathbb{K}^{k} \left(A,U,r_0 \right).
											\end{align*}
											The final solution update and the residual recurrences take the following form:
											\begin{align*}
												x_{k+1} = x_k + \alpha_{k} p_k, \\
												r_{k+1} = r_k + \alpha_k A p_k,
											\end{align*}	
											where 
											\begin{align*}
												p_{k} & = \beta_{k-1} p_{k-1} + \left(I-  U(U^TAU)^{-1} (AU)^T \right) r_{k}, \\
												\alpha_{k} & = \left(r_{k}^T r_{k}\right)/ \left(p_{k}^T A p_{k}\right), \\
												\beta_{k-1} & = \left(r_{k}^T r_{k}\right) / \left(r_{k-1}^T r_{k-1}\right).
											\end{align*}
											
											Next, we discuss how to use the above machinery	for our requirements. Consider solving the linear system given by (\ref{eq:inexcat2a}), originally (\ref{inexcat2a}). For adapting the Lanczos process in Section \ref{sec:Lanczos_exta_orth}, while solving this linear system, we need to achieve the extra orthogonality in (\ref{eq:ext_Lan_1st}). Similarly, for adapting  the Ritz-Galerkin projection in Section \ref{sec:Ritz-galerkin_extra_orth}, while solving this linear system, we need to achieve the extra orthogonality in  (\ref{eq:second_resid_AIRGA1}). Both these orthogonalities can be achieved if we take
											{\begin{align}\label{eq:1st_recycle_col}
													U = \left[\eta_{(0)} \ \tilde{V}_1\right]
												\end{align}}
												in RCG.
												
												By defining $U$ as above, {$\eta_{(0)}$} and $\tilde{V}_1$ are added in the Krylov search space, which is not needed in the adapted Lanczos process.
												Also, we are doing extra work here since {$\eta_{(0)}$} orthogonality is needed only for Lanczos (not for Ritz-Galerkin), and $\tilde{V}_1$ is needed for Ritz-Galerkin (not for Lanczos).  
												
												These facts are true but besides the benefit of ease of implementation, this choice of space often leads to acceleration of the system. We support this with experiments in the next section. A theoretical study of this choice of space is the part of future work.
												
												{Also, to satisfy all the other orthogonalities of the previous subsection, equivalent of $U$  (as in \eqref{eq:1st_recycle_col}) can be defined}.
												{Since we are usually more concerned about the accuracy of the obtained reduced dynamical systems, we investigate this apects next.}
												
												
												
												
												
												
												\section{Accuracy of the Reduced Systems}\label{sec:accuracy}
												Using Theorem 15.1 of \cite{trefethen1997numerical} 
												we know that if the AIRGA algorithm is backward stable, then the relative accuracy of the reduced system obtained by using the inexact AIRGA algorithm, as compared to using the exact AIRGA algorithm, is given as follows:
												\begin{align}
													\label{eq:back_stab}
													\frac{\|\hat{H}(s)- \tilde{\hat{H}}(s)\|_{H_{2}}}{\|\hat{H}(s)\|_{H_{2}}} = \mathcal{O}\left(\kappa \left(H(s)\right)\cdot\|Z\|_2 \right),
												\end{align}
												where $\kappa \left(H(s)\right)$ is the condition number of $H(s)$ (discussed below), and $Z$ is the perturbation in $H(s)$. 
												As earlier, $\hat{H}(s) $ is the reduced system obtained by using the exact AIRGA algorithm and $\tilde{\hat{H}}(s)$ is the reduced system obtained by using the inexact AIRGA algorithm. We are looking at reduced systems obtained at line \ref{AIRGA:red_sys}  of Algorithm~\ref{Algo:AIRGA}. 
												That is, after each step of the outer \texttt{while} loop (line \ref{AIRGA:while-outer}). 
												Thus, accuracy of the reduced system is dependent on the conditioning of the problem as well as the perturbation. Next, we look at both these quantities separately.
												
												\subsection{Conditioning Expression}
												We want to compute conditioning of our system with respect to performing the inexact linear solves on lines \ref{algo:line-1st} and \ref{algo:line-jth} of Algorithm~\ref{Algo:AIRGA}. Since for backward stability we equate the reduced system obtained by performing the inexact AIRGA algorithm on the unperturbed (original) full system $\left(H(s)\right)$ and performing the exact AIRGA algorithm on the perturbed full system  $\left(\tilde{H}(s)\right)$, these inexact linear solves are captured by $\tilde{H}(s)$. Thus, the conditioning of the input dynamical system with respect to computing the $H_2$-norm of the error system $H(s) - \tilde{H}(s)$ will give us a \textit{good approximation} to the 
												conditioning of the input dynamical system that we want to access
												$\left( \text{with respect to computing the} \ H_2-\text{norm of} \ \hat{H}(s)-\tilde{\hat{H}}(s) \right)$. Similar behaviour has been observed for linear first-order   dynamical systems $\left(\text{see Theorem 3.1 and 3.3 in \cite{Beattie20122916}} \right)$ and  bilinear first-order  dynamical systems \cite{CHOUDHARY201856}. 
												
												Recall, the condition number by definition means the relative change in the output  $\left(\text{for us this is} \ \|H(s) - \tilde{H}(s)\|_{H_2}/\|H(s)\|_{H_2} \right)$ with respect to the relative change in the input $\left( \text{for us this is} \ \|Z\|_2/\|K\|_2 \ \text{since we are perturbing the} \ K \ \text{matrix} \right)$ \cite{CHOUDHARY201856}. Hence, from (\ref{eq:sec_con_bound})  we have
												\begin{align}
													\frac{\|H(s)-\tilde{H}(s)\|_{H_2}}{\|H(s)\|_{H_2}} \le \frac{\|A(s)^{-1}B\|_{H_{\infty}}\|{\mathscr C}^TA(s)^{-1}\|_{H_{2}}}{\|H(s)\|_{H_2}} \cdot \frac{\|K\|_2}{1-\|A(s)^{-1}\|_{H_\infty}} \cdot \frac{\|Z\|_2}{\|K\|_2},
												\end{align}
												where it is assumed that $\|Z\|_2<1$ and $\|A(s)^{-1}\|_{H_\infty} <1$. Hence, the above inequality is equivalent to
												\begin{align}
													\frac{\|H(s)-\tilde{H}(s)\|_{H_2}}{\|H(s)\|_{H_2}} \le \kappa\left(H(s)\right) \cdot \frac{\|Z\|_2}{\|K\|_2},
												\end{align} 
												where, \begin{align}\label{eq:cond_num_prob}
													\kappa \left(H(s)\right) = \frac{\|A(s)^{-1}B\|_{H_{\infty}} \|{\mathscr C}^TA(s)^{-1}\|_{H_{2}}}{\|H(s)\|_{H_2}} \cdot \frac{\|K\|_2}{1-\|A(s)^{-1}\|_{H_\infty}}. 
												\end{align}
												
												In the numerical experiments section (Section \ref{sec:Numerical_results}), for the first example taken, we show that this condition number is fairly small, whereas, for the second one it is large. In other words, the first problem is well conditioned and the second problem  is ill-conditioned
												with respect to the $H_2$-norm of the error system $H(s)-\tilde{H}(s)$\footnote{This ill-conditioning of the second problem does not effect our main conjecture. We discuss this aspect in-detail later. }. Note that $\|Z\|_2 < 1$ and $\|A(s)^{-1}\|_{H_\infty} < 1$, as assumed here, come from the assumptions for backward stability of the AIRGA algorithm (see Theorem 2), and hence, we do not need any extra assumptions.
												
												{
													\subsection{Computation of Perturbation}
													Recall (\ref{eq:Orth_X}), which has the form
													\begin{align}\label{eq:prertub_Z}
														Z \mathbf{X} = \eta.	
													\end{align}
													Here, $Z \in \mathbb{R}^{n \times n}, \mathbf{X} \in \mathbb{R}^{n \times mJ}$, and $\eta \in \mathbb{R}^{n \times mJ}$. Also note that we are solving for $Z$. As discussed in Introduction, 
													the upper bound for $J$ is $\lceil r_{\max}/m \rceil$, and hence, $mJ \leq r_{\max}$. Using the fact that $r_{\max} \ll n$, we have $mJ < n$. Hence, we have an underdetermined system of equations, which will have more than one solution. For such a system, Singular Value Decomposition (SVD) helps provide one solution \cite{Golub1996}. This SVD for $\mathbf{X}$ is given as follows:
													\begin{align*}
														\mathbf{X} = \mathbb{U} \Sigma \mathbb{V}^T,
													\end{align*}
													where $\mathbb{U} \in \mathbb{R}^{n \times n}, \mathbb{V} \in \mathbb{R}^{mJ \times mJ}$ are unitary matrices (i.e. $\mathbb{U} \mathbb{U}^T =I$ and  $\mathbb{V}\mathbb{V}^T =I $) and $\Sigma \in \mathbb{R}^{n \times mJ}$ is a diagonal matrix comprising of singular values of $\mathbf{X}$. Let ${r_{n}} = rank(\mathbf{X})$, then $\Sigma = \diag \left(\varsigma_1, \ldots, \varsigma_{r_n},0, \ldots, 0 \right)$. Partitioning $\mathbb{U}$ as $ [\mathbb{U}_1 \ \mathbb{U}_2]$ and  $\mathbb{V}$ as $[\mathbb{V}_1 \ \mathbb{V}_2]$, where $\mathbb{U}_1$, $\mathbb{V}_1$ have $r_n$ columns; $U_2$, $V_2$ have the remaining columns of $U$, $V$, respectively; and  $r_n \le mJ$, we get
													\begin{align}\label{eq:svd_X1}
														\mathbf{X} = \begin{bmatrix}
															\mathbb{U}_1 & \mathbb{U}_2
														\end{bmatrix}\begin{bmatrix}
														\Sigma_{r_n} & 0 \\
														0 & 0
													\end{bmatrix} \begin{bmatrix}
													\mathbb{V}_1 & \mathbb{V}_2
												\end{bmatrix}^T,
											\end{align}
											where $\Sigma_{r_n} = \diag\left(\varsigma_1, \ldots \varsigma_{r_n}\right)$. 
											By using (\ref{eq:svd_X1}) and definition of Moore-Penrose Pseudoinverse (\cite{Meyer2000}; page 423) we have
											\begin{align*}
												\mathbf{X}^{\dagger} = \begin{bmatrix}
													\mathbb{V}_1 & \mathbb{V}_2
												\end{bmatrix} \begin{bmatrix}
												\Sigma_{r_n}^{-1} & 0 \\
												0 & 0
											\end{bmatrix}\begin{bmatrix}
											\mathbb{U}_1^T \\ \mathbb{U}_2^T
										\end{bmatrix}.
									\end{align*} 
									Substituting the above expression in  (\ref{eq:prertub_Z}), we have\footnote{{If the system $Z {\mathbf{X}} = {\eta}$ is in-consistent, then this is the least squares solution.}} 
									\begin{align}\label{eq:Z_svd}
										Z = {\eta} \mathbb{V}_1 \Sigma^{-1}_{r_n} \mathbb{U}_1^T.
									\end{align}
								}
								
								Next, we relate the perturbation $Z$ and the  cumulative residual ${\eta}$.
								{\begin{align}\label{delta_K_norm1}
										\|Z\|_2 \le \|Z\|_{f} &\le \|{\eta} \cdot \mathbb{V}_1 \Sigma_{r_n}^{-1} \mathbb{U}_1^T\|_{f} 
										\le\|{\eta}\|_{f} \|\mathbb{V}_1 \Sigma_{r_n}^{-1} \mathbb{U}_1^T\|_{f}, \\ \notag
										& \le \left(\|-{\eta}_{(0)}\|_f + \cdots + \|-{\eta}_{(J-1)}\|_f\right) \left(\|\mathbb{V}_1 \Sigma_{r_n}^{-1} \mathbb{U}_1^T\|_{f} \right).
									\end{align}}
									In the above equation, {$-{\eta}_{{(0)}}, \ldots, -{\eta}_{(J-1)}$} 
									represent the residuals obtained while solving the linear systems arising in the model reduction process. These residuals will reduce if we solve such linear systems more accurately.
									The second term  $\|\mathbb{V}_1 \Sigma_{r_n}^{-1} \mathbb{U}_1^T\|_{f}$ is usually more dependent on the selection of the expansion points ($s_{i}$), and less on the accuracy to which we solve the linear systems~\cite{Beattie20122916}. We support this argument with numerical experiments. 
									
									
									To summarize from \eqref{eq:back_stab} we know, $\|\hat{H} (s)-\tilde{\hat{H}}(s)\|_{H_2}$ is proportional to $\kappa \left(H(s)\right)$ and $\|Z\|_2$. The problem is usually well conditioned and $\|Z\|_2$ is directly proportional to the cumulative residual norm $\|\eta\|_f$ (as in \eqref{delta_K_norm1}). Thus, assuming backward stability conditions hold (as discussed in the previous section), as we iteratively solve the linear systems arising in the AIRGA algorithm more accurately (i.e. reduce the stopping tolerance of the linear solver), we should get a more accurate reduced system. This is very useful in deciding  when to stop the linear solver. If we need a very accurate reduced system, then we need to iterate more in the linear solver, else we can stop earlier. 
									

\section{Numerical Experiments}
\label{sec:Numerical_results}
		
{
	As motivated in Section~\ref{sec:Sat_Back_stab_cond}, for stability we focus on the CG method for solving the linear systems arising  in the AIRGA algorithm. Also, as discussed earlier,  CG is optimal for SPD linear systems. Thus, we need to ensure that the coefficient matrices of all the linear systems to be solved are SPD. 
	
	The coefficient matrices are of the form $s^2_i M + s_i D+ K$. To achieve that these matrices are SPD at start we do as below.
	
	(a) We take input models that have M, D and K matrices as SPD. We use the one dimensional beam model (size $10,000$) \cite{BeattieG2005} and the Gyroscope model (size $17,361$) \cite{BGyroscopemodel} that have such matrices and are commonly used (discussed in the next two subsections). These models are of the form~\cite{BeattieG2005,BGyroscopemodel}
	\begin{align}\label{eq:sec_dy_sys}
		\begin{split}
			& M \ddot{x}(t)+D\dot{x}(t)+Kx(t) = Bu(t), \\ 	
			& y(t)=Cx(t),		
		\end{split}
	\end{align}
	where 
	$M, \ D, \  K  \in  \mathbb{R}^{n \times n} $ are the mass, the damping and the stiffness matrices, respectively,  
	$B  \in  \mathbb{R}^{n \times 1} \ \text{and} \ C   \in  \mathbb{R}^{1 \times n}$. These models are Single Input Single Output (SISO), and have proportional damping, i.e. $D = \alpha M+\beta K$, where the damping coefficients $\alpha$ and $\beta$ belong to $(0,1)$.

	(b) We take the input expansion points $(s_i)$ to be real and positive. 
	
	Next, we discuss how to ensure that the linear system matrices are SPD after the first AIRGA iteration (i.e. after start).
	After the first AIRGA iteration, the expansion points are chosen from the eigenvalues of the quadratic eigenvalue problems of the form  $\lambda^2 \hat{M} + \lambda \hat{D} + \hat{K}$. For  both our models, these eigenvalues  turn out to be complex (case $3.8$ of Table 1.1 in~\cite{Tisseur2001}). 
	Thus, we get complex expansion points. 
	Execution of the AIRGA algorithm as well as the accuracy of the reduced system does not get affected if one uses real expansion points or complex expansion points. 
	Since real expansion points here are positive too (again because of case $3.8$ of Table 1.1 in~\cite{Tisseur2001}), using them ensures that our coefficient matrices, $s_i^2 M +s_i D+ K$, are SPD at all the AIRGA iterations. Hence, we use real expansion points.
	
	In Algorithm~\ref{Algo:AIRGA}, at line \ref{AIRGA:while-outer}, the overall iteration (\texttt{while-loop}) terminates when the change in the reduced system (computed as the $H_{2}$-error between the reduced systems of two consecutive AIRGA iterations) is less than a certain tolerance. We take this tolerance to be $10^{-04}$ based on values in~\cite{Bonin20161}. There is one more stopping criteria in this algorithm at line \ref{AIRGA:while-inner}.  This checks the $H_{2}$-error between two temporary reduced systems. We take this tolerance to be $10^{-06}$ based upon values in~\cite{Bonin20161}.
	
	As also motivated in Section \ref{sec:Sat_Back_stab_cond}, to ensure that
	the extra orthogonalities for a stable AIRGA algorithm are satisfied, we use RCG instead of CG. As earlier, we refer to this as the inexact AIRGA algorithm. 
	Preconditioning has to be employed when iterative methods fail or have a very slow convergence. Here, for the first model, we observe that the unpreconditioned RCG method has slow convergence whereas in the second model it fails to converge. Thus, we use a preconditioner. Since Sparse Approximate Inverse (SPAI) \cite{Chow1998}  and  Incomplete Cholesky Factorization (ICHOL) \cite{Manteuffel1980, Saad2003} are the most general types of preconditioners, 
	we can use any of these preconditioners with RCG. 
	Here, we use the standard SPAI (with stopping tolerance of $10^{-04}$) for the first model and the standard ICHOL (with drop tolerance of $10^{-04}$) for the second model.
	For comparison, we solve all linear systems by a direct method as well. As earlier, we refer to this as the exact AIRGA algorithm. For certain types of analyses, we compare CG and RCG behaviours too.
	
	We implement our codes in MATLAB (2016b), and test on a machine with the following configuration: Intel Xeon(R) CPU~E5-1620~V3~@~3.50 GHz., frequency 1200 MHz., 8 CPU and 64 GB RAM. 
}

\subsection{One Dimensional Beam Model} \label{1dbeam}
{
	As discussed earlier, we do experiments on a system of size $10,000$. Damping coefficients $\alpha$ and $\beta$ both are taken as $0.05$ \cite{BeattieG2005} {we take three expansion points as $s_1 = 0.3142, s_2 = 0.6283$ and $s_3 = 0.9425$ (based upon experience)}. The maximum dimension to which we want to reduce the system ($r_{max}$) is taken as {$3$ based upon experience}. 
	Thus, in the AIRGA algorithm, we have to solve linear systems of size $10,000 \times 10,000$. While using RCG for solving these linear systems, we use two different stopping tolerances {$10^{-02}$ and $10^{-08}$}. Ideally, as discussed earlier, we should obtain a more accurate reduced system for the smaller stopping tolerance.

	First, let us look at the assumptions for backward stability of the AIRGA algorithm (see Theorem 2). 
	{While referring to this theorem, we have already satisfied the conditions (a) and (b) by using CG and RCG, respectively. At all AIRGA iterations, $\sigma_1$ picked is $s_1$, $\sigma_2$ picked is $s_2$, and $\sigma_3$ picked is $s_3$. Thus, at all AIRGA iterations, we solve the linear systems in the following order to match our theory proposed in subsections \ref{sec:Lanczos_exta_orth} and \ref{sec:Ritz-galerkin_extra_orth}:
		\begin{itemize}
			\item  \eqref{inexact1} (including \eqref{eq:inexact1})-- Three systems corresponding to three expansion points
			\item \eqref{inexcat2a} (or \eqref{eq:inexcat2a})-- One system
			\item \eqref{inexact1} (or \eqref{eq:inexact1_rep})-- One system {\it resolve}
			\item \eqref{inexact2} (or \eqref{eq:inexact3a})-- One system
			\item \eqref{inexact1} (or \eqref{inexact1-resol})-- One system {\it resolve}
		\end{itemize}
	}
	
	{Hence, at all AIRGA iterations instead of solving five linear systems, we solve seven linear systems. This is acceptable because this gives us a stable MOR algorithm. Sometimes, use of recycle space accelerates the convergence of all linear systems in-turn offsetting this extra cost. We demonstrate this behaviour in the next example.}
	
	{Next, we analyze the assumptions (c) and (d) of Theorem 2}.
	For all expansion points, $A(s)$ is invertible and $\|A(s)^{-1}\|_{H_\infty}$ is less than one. E.g., for the initial set of expansion points, $\|A(s)^{-1}\|_{H_\infty}$ is $2.68 \times 10^{-02}$. Finally, $\|Z\|_2$, at the end of the first AIRGA iteration, for the RCG stopping tolerance of {$10^{-02}$ and $10^{-08}$ is $0.28$ and $ 0.06$,} respectively, both of which are also less than one \footnote{{Our ${\mathbf{X}}$ while solving (\ref{eq:Orth_X}) using (\ref{eq:Z_svd}) is full rank i.e. 8 here.}}.
	These values are  less than one at the end of all the other AIRGA iterations as well. The condition number for our problem, as defined in (\ref{eq:cond_num_prob}), is $8.63 \times 10^{-02}$. This shows that the one dimensional beam model is well-conditioned. {As earlier, we still use the SPAI preconditioner for better acceleration.}

	The accuracy results are given in Fig. \ref{fig:acc_ana_beam}. 
	{We use the following settings: expansion points $s_i = 2 \pi f$, where frequency $f$ vector consists of equally spaced twenty points between $25$ and $250$.}
	In Fig. \ref{fig:acc_ana_beam}, we have the accuracy of the reduced system $\left(\|\hat{H}(s)-\tilde{\hat{H}}(s)\|_{H_2}\right)$ on the y-axis and the AIRGA iterations on the x-axis.
	Here, the dotted line corresponds to the RCG stopping tolerance of {$10^{-02}$} while the solid line corresponds to the RCG stopping tolerance of {$10^{-08}$}.
	From Fig. \ref{fig:acc_ana_beam}, it is evident that we get a more accurate reduced system as we solve the linear systems more accurately (dotted line is above the solid one at all the AIRGA iterations). 
	
	{In Table \ref{tab:accuracy_anl}, we give the accuracy results corresponding to each AIRGA iterations. The AIRGA algorithm gets more consistent as it converges to ideal expansion points. Hence, the accuracy of the reduced system for the RCG stopping tolerance of $10^{-08}$ is visibly better than the accuracy of the reduced system for the RCG stopping tolerance of $10^{-02}$.
	}
	
	
	\begin{figure}[]
		\includegraphics[width=140mm]{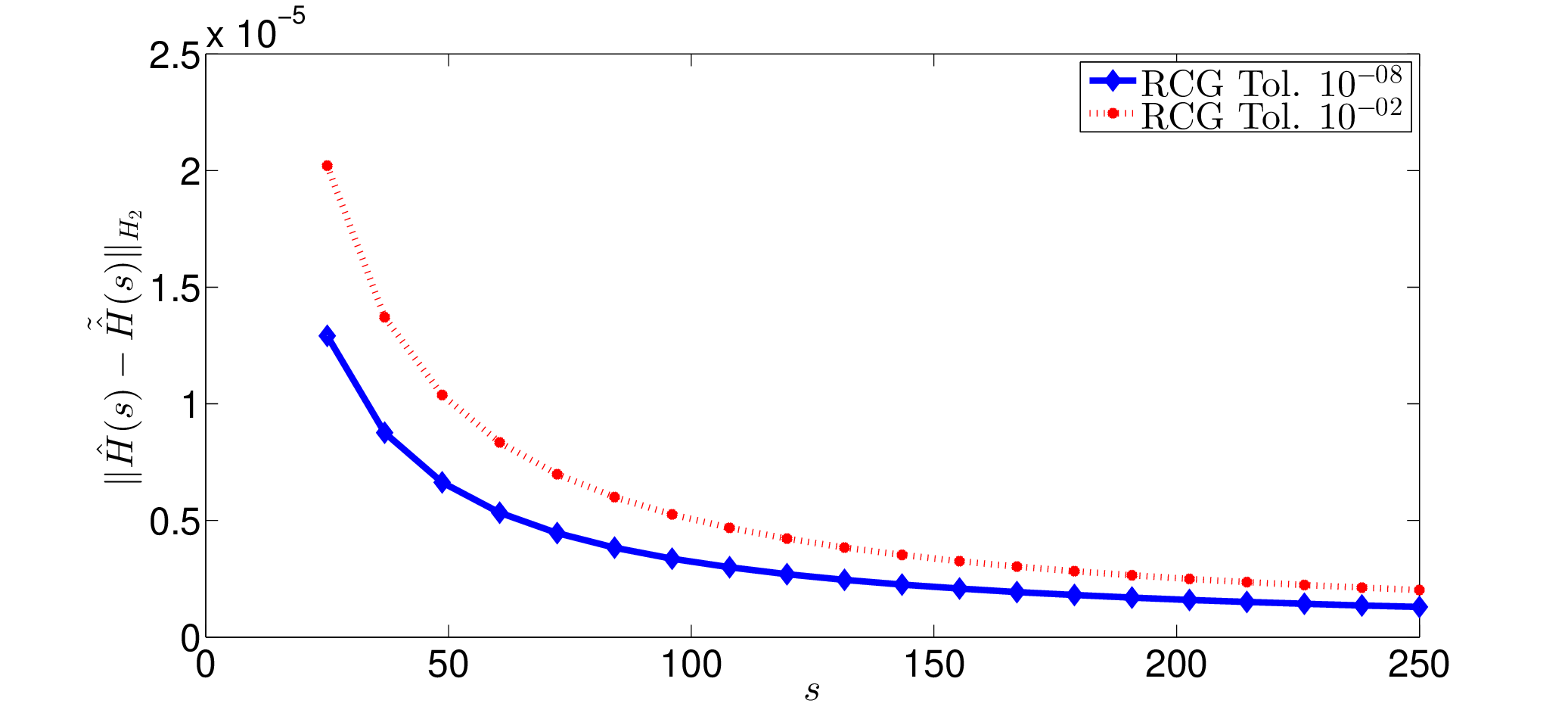}
		\caption{Accuracy of the reduced system plotted at each AIRGA iteration for two different stopping tolerances in RCG; one dimensional beam model.}
		\label{fig:acc_ana_beam}       
	\end{figure}
	\begin{table}[]
		\centering
		\caption{Accuracy of the reduced system at each AIRGA iteration for the two different stopping tolerances in RCG; one dimensional beam model.}
		\label{tab:accuracy_anl}
		\begin{tabular}{|c|c|c|}
			\hline
			\multirow{2}{*}{\begin{tabular}[c]{@{}c@{}}AIRGA \\ Iteration\end{tabular}} & \multicolumn{2}{c|}{$||\hat{H} - \tilde{\hat{H}}||_{H_{2}}$} \\ \cline{2-3} 
			& RCG stopping tolerance of $10^{-02}$                    &  RCG stopping tolerance of $10^{-08}$  \\ \hline
			1 & $8.34 \times 10^{-05}$         & $2.99 \times 10^{-06}$        \\ \hline
			2 & $6.98 \times 10^{-05}$         & $2.69 \times 10^{-06}$        \\ \hline
			3 & $1.37 \times 10^{-05}$         & $2.45 \times 10^{-06}$        \\ \hline
			4 & $1.03 \times 10^{-05}$         & $2.25 \times 10^{-06}$        \\ \hline
		\end{tabular}
	\end{table}

	\subsection{Gyroscope Model}\label{Gyroscope}
	{
		As mentioned earlier, we do another experiment on a system of size $17,361$. Damping coefficients $\alpha$ and $\beta$ are taken as $0.2$ and $1.34 \times 10^{-04}$, respectively \cite{BGyroscopemodel}. Here, again, {we take three expansion points as $s_1 = 6.2832, s_2 = 317.3009$ and $s_3 = 628.3185$.}
		The  dimension to which we want to reduce the system ($r_{max}$) is taken as $12$ based upon similar values in \cite{BGyroscopemodel}. Here, in the AIRGA algorithm we have to solve the linear systems of size $17,361 \times 17,361$. Again, we use RCG for solving these linear systems. {To demonstrate our main result, we ideally want the stopping tolerances to be six orders of magnitude different from each other. E.g., $10^{-02}$ and $10^{-08}$ in the previous problem. Here, we are unable to solve the linear systems for tolerances less than $10^{-10}$. As for the higher tolerance, if we go beyond $10^{-08}$, then the AIRGA algorithm's convergence varies (differing iteration counts for convergence). Thus, we cannot compare results of the two cases. Hence, we use stopping tolerances of $10^{-08}$ and $10^{-10}$.}
		As discussed earlier, we should obtain a more accurate reduced system for the smaller stopping tolerance.    
		
		Similar to the previous experiment, here also we look at the remaining assumptions for  backward stability of the AIRGA algorithm (see Theorem 2). 
		{While referring to this theorem, we have already satisfied the conditions (a) and (b) by using CG and RCG, respectively. When applying theorem of Section \ref{sec:back_gen_con} here (to satisfy (b)), for simplicity, we do not perform the required resolves. The results below demonstrates that this approximation does not have any effect on our intended behaviour. We still get a more accurate reduced system as we solve the linear systems more accurately. Here, we also show that use of recycle space accelerate the convergence of the linear systems.   
		}
		
		{Next, again, we analyze the assumptions (c) and (d) of Theorem 2}.	
		For all expansion points, $A(s)$ is invertible and $\|A(s)^{-1}\|_{H_\infty}$ is less than one. E.g., for the initial set of expansion points, $\|A(s)^{-1}\|_{H_\infty}$ is $6.46 \times 10^{-01}$. Finally, $\|Z\|_2$, at the end of the first AIRGA iteration, for the RCG stopping tolerance of $10^{-08}$ and $10^{-10}$ is $8.6 \times 10^{-01}$ and $3.3 \times 10^{-01}$, respectively, both of which are also less than one \footnote{{ Our ${\mathbf{X}}$ while solving (\ref{eq:Orth_X}) using (\ref{eq:Z_svd}) is rank deficient (10 instead of 12) but that does not affect our computations.}}. 
		These values are  less than one at the end of all the other AIRGA iterations as well. 
		
		The condition number for this problem, as defined in (\ref{eq:cond_num_prob}), is $5.15 \times 10^{09}$. 
		This shows that the Gyroscope model is ill-conditioned. As earlier, we use the basic ICHOL preconditioner, which helps to reduce the amount of ill-conditioning but does not completely eliminated it. If needed, we can use a more advanced preconditioner.

		Accuracy of the reduced system is proportional to the condition number $\kappa\left(H(s)\right)$ and the perturbation $\|Z\|$ (see \eqref{eq:back_stab}). Since, the condition number here {still remain} high, we  get a less accurate reduced system. However, this is still a good problem for us since we want to demonstrate that the reduction in perturbation (linked to linear solver stopping tolerance) improves accuracy. High condition number spoils the accuracy equally for both the RCG stopping tolerances ($10^{-08}$ and  $10^{-10}$). The accuracy results are given in  Table \ref{Tab:accur_anal_17_1K}. It is again evident that we get a more accurate reduced system as we solve the linear systems more accurately. 
		


		For this model, we observe that the number of iterations required for convergence of RCG is less than that of CG, both of which are given in Table \ref{Tab:cong_cg_RCG_anal_17K}.    
		We see a savings of about $10 \%$ in the average linear solver iterations. The corresponding computation times are given in Table \ref{Tab:cong_cg_RCG_anal_time_17K}. The savings in iterations translate to about $5 \%$ savings in time.

		Here, we do some other analysis corresponding to (\ref{delta_K_norm1}), i.e. relation between the perturbation and the stopping tolerance. From  Table \ref{Tab:pertu_anal_17_1K}, we demonstrate that $\|\mathbb{V}_1 \Sigma_{r_n}^{-1} \mathbb{U}_1^T\|_f$ is less sensitive to the accuracy to which we solve the linear systems (as we reduce the stopping tolerance of RCG from $10^{-08}$ to $10^{-10}$, $\|\mathbb{V}_1 \Sigma_{r_n}^{-1} \mathbb{U}_1^T\|_f$ stays almost the same). 
	}
	\begin{table}[]
		\centering
		\caption{Accuracy of the reduced system at each AIRGA iteration for the two different stopping tolerances in RCG; Gyroscope Model.}
		\label{Tab:accur_anal_17_1K}
		\begin{tabular}{|c|c|c|}
			\hline
			\multirow{2}{*}{\begin{tabular}[c]{@{}c@{}}AIRGA \\ Iteration \end{tabular}} & \multicolumn{2}{c|}{$||\hat{H}-\tilde{\hat{H}}||_{H_2}$} \\ \cline{2-3} 
			& RCG stopping tolerance $10^{-08}$ & RCG stopping tolerance $10^{-10}$ \\ \hline
			1 & $ 1.55 \times 10^{-03}$ & $8.66 \times 10^{-04}$  \\ \hline
			2 & $3.63 \times 10^{-05}$ &  $3.14 \times 10^{-05}$ \\ \hline
		\end{tabular}
	\end{table}
	
	
	%
	%
	
	\begin{table}[]
		\centering
		\caption{Convergence analysis of CG and RCG at two different stopping tolerances; Gyroscope Model.}
		\label{Tab:cong_cg_RCG_anal_17K}
		\begin{tabular}{|c|c|c|c|c|}
			\hline
			\multirow{2}{*}{\begin{tabular}[c]{@{}c@{}}AIRGA \\ Iteration \end{tabular}} & \multicolumn{2}{c|}{Stopping tolerance $10^{-08}$} & \multicolumn{2}{c|}{Stopping tolerance $10^{-10}$} \\ \cline{2-5} 
			& Avg. CG Itr. &  Avg. RCG Itr. & Avg. CG Itr. & Avg. RCG Itr.	\\ \hline   
			1 & 216 & 207 & 244 & 224  \\ \hline
			2 & 202 & 180 & 228 & 206 \\ \hline		
			\textbf{Total} & \textbf{418} & \textbf{387} & \textbf{472} & \textbf{430}  \\ \hline
		\end{tabular}
	\end{table}
	
	\begin{table}[]
		\centering
		\caption{Computation time of CG and RCG at two different stopping tolerances; Gyroscope Model.}
		\label{Tab:cong_cg_RCG_anal_time_17K}
		\begin{tabular}{|c|c|c|c|c|}
			\hline
			\multirow{2}{*}{AIRGA} & \multicolumn{2}{c|}{\begin{tabular}[c]{@{}c@{}}Stopping tolerance $10^{-08}$  \end{tabular}}  & \multicolumn{2}{c|}{\begin{tabular}[c]{@{}c@{}}Stopping tolerance $10^{-10}$  \end{tabular}} \\  \cline{2-5} 
			Iteration & {\begin{tabular}[c]{@{}c@{}} CG time \\ (secs.) \end{tabular}} & {\begin{tabular}[c]{@{}c@{}} RCG time \\ (secs.) \end{tabular}} &  {\begin{tabular}[c]{@{}c@{}} CG time \\ (secs.) \end{tabular}} & {\begin{tabular}[c]{@{}c@{}} RCG time \\ (secs.) \end{tabular}} \\ \cline{1-5}  		
			1 & 2.35 & 2.20 & 2.49 &   2.41 \\ \hline
			2 & 2.04 & 1.95 & 2.33 &   2.23 \\ \hline
			\textbf{Total} & \textbf{4.39} & \textbf{4.15} & \textbf{4.82} &   \textbf{4.64} \\ \hline
		\end{tabular}
	\end{table}
	
	\begin{table}[]
		\centering
		\caption{The perturbation expression quantities for RCG at two different stopping tolerances; Gyroscope Model.}
		\label{Tab:pertu_anal_17_1K}
		\begin{tabular}{|c|c|c|c|c|}
			\hline
			\multirow{2}{*}{\begin{tabular}[c]{@{}c@{}}AIRGA \\ Iteration \end{tabular}} & \multicolumn{2}{l|}{RCG Stopping tolerance $10^{-08}$}  & \multicolumn{2}{l|}{RCG Stopping tolerance $10^{-10}$} \\ \cline{2-5} 
			& $\|\eta\|_f$ & $\|\mathbb{V}_1 \Sigma_{r_n}^{-1} \mathbb{U}_1^T\|_f$ & $\|\eta\|_f$  & $\|\mathbb{V}_1 \Sigma_{r_n}^{-1} \mathbb{U}_1^T\|_f$ 	\\ \hline
			1 & $2.5 \times 10^{-09}$ & $1.34 \times 10^{09} $ & $2.9 \times 10^{-10}$ & $1.33 \times 10^{09}$   \\ \hline
			2 & $2.6 \times 10^{-09} $ & $1.98 \times 10^{11} $ & $2.6 \times 10^{-10}$ & $1.98 \times 10^{11}$ \\ \hline
		\end{tabular}
	\end{table}
	
\section{Conclusion and Future Works}
\label{sec:Conc_future}

{
We discuss application of preconditioned iterative methods for solving the large linear systems {in the AIRGA algorithm. This algorithm  is used for reducing linear non-parametric second-order dynamical systems with proportional damping. These iterative} methods find solutions only up to a certain tolerance. {Hence, we show that} under {four} mild conditions, AIRGA is backward stable with respect to these inexact linear solves. We also analyze the accuracy of the resulting reduced system, and support all our results with multiple numerical experiments.





The \textit{first condition} for stability enforces the use of a Ritz-Galerkin based linear solver, where the
residual of a linear system is made orthogonal to the corresponding Krylov subspace. The \textit{second condition} for stability requires satisfying few other orthogonalities. Since the CG method is the most popular linear solver based upon the Ritz-Galerkin theory and is ideal for SPD linear systems, we focus on SPD systems only. We use Recycling CG (RCG) to achieve the extra orthogonalities. The future work here involves modifying other methods based upon the Ritz-Galerkin theory (to achieve extra orthogonalities), which can be used to solve general non-symmetric indefinite linear systems. For example, the Full Orthogonalization Method (FOM).

The \textit{third condition} for stability {involves computing} $A(s)$, which is a function of the frequency $(s)$ and {input} dynamical system matrices, 
{inverting it, and bounding its norm by one}.
{This condition is easily} satisfied for all our models, but {it} may not always hold. Future work here involves better characterizing {this condition} in-terms of the underlying dynamical system.

The \textit{fourth and final condition} for stability involves being able to compute perturbation $Z$ from the given expression and bounding its norm by one. As earlier, although for all our models {this condition is} easily satisfied, {it} may not always hold. $Z$ is dependent on the linear solver stopping tolerances. Hence, we need to study range of these tolerances when the norm of this perturbation could be bounded by one.

The condition number of the dynamical system, which we use is an approximation to the ideal condition number. That is, condition number of the dynamical system
with respect to computing the $H_2$-norm of the error between the inexactly computed reduced system and the exactly computed reduced system. This is also part of the
future work.	
}
\ \\ 
\ \\

\section*{Acknowledgments}
We would like to deeply thank Prof. Heike Fassbender (at TU Braunschweig, Germany) for many fruitful discussions on this project. 

\bibliographystyle{siamplain}
\bibliography{StabilityAIRGA}

\end{document}

%% file: ex_shared.tex
\usepackage{lipsum}
\usepackage{amsfonts}
\usepackage{graphicx}
\usepackage{epstopdf}
\usepackage{caption}
\usepackage{multirow}
\usepackage{array, multirow}
\usepackage{threeparttable}
\usepackage{xcolor,colortbl}
\usepackage{amsmath,amssymb}
\usepackage{amsmath,amssymb,amsfonts,epsfig,color,graphicx,eucal,mathrsfs}
\usepackage{amsmath}
\usepackage{epstopdf} 
\usepackage[utf8]{inputenc}
\usepackage[multiple]{footmisc}
 \usepackage{upgreek}
\usepackage{lipsum}
\usepackage{xcolor}

\RequirePackage{algorithm}
\usepackage{algcompatible}

\definecolor{mygreen}{RGB}{10,150,10}
\definecolor{myblue}{RGB}{10,10,150}

\ifpdf
  \DeclareGraphicsExtensions{.eps,.pdf,.png,.jpg}
\else
  \DeclareGraphicsExtensions{.eps}
\fi

\numberwithin{theorem}{section}

\newcommand{\TheTitle}{Stability Analysis of Inexact Solves in Model Reduction} 
\newcommand{\TheAuthors}{Kapil Ahuja and Navneet Pratap Singh}

\headers{\TheTitle}{\TheAuthors}

\title{Stability Analysis of Inexact Solves in Model Reduction of Non-parametric Second-order Dynamical systems}

\author{
Kapil Ahuja\thanks{Math of Data Science and Simulation (MODSS) Lab, Department of Computer Science and Engineering, Indian Institute of Technology Indore, India (\email{kahuja@iiti.ac.in})}
\and
Navneet Pratap Singh\thanks{Bennett University, Greater Noida, Uttar Pradesh, India (\email{Navneet.Singh@bennett.edu.in)}}
}

\usepackage{amsopn}
\DeclareMathOperator{\diag}{diag}
